\documentclass[journal]{IEEEtran}
\newlength{\imgheight}
\newlength{\imgheighth}
\usepackage{amsmath,amssymb,amsfonts}
\usepackage[linesnumbered,ruled,vlined]{algorithm2e}
\usepackage{graphicx}
\graphicspath{ {./figs/} }
\usepackage{multirow}
\usepackage{booktabs}
\usepackage{subfigure}
\usepackage{stfloats}
\usepackage{array}
\usepackage{tabularx}
\usepackage{hhline}
\usepackage{xspace}
\usepackage{upgreek}
\usepackage{color}
\usepackage{makecell}
\usepackage{mathtools}
\usepackage{xcolor}
\usepackage{spverbatim}

\newcounter{itemlistc}

\newenvironment{itemlist}
    {   \setcounter{itemlistc}{0}
    \begin{list}{$\bullet$}
        {\usecounter{itemlistc}
        \setlength{\parsep}{0pt}
        \setlength{\topsep}{3pt}
        \setlength{\itemsep}{0pt}}
        }{ \end{list} }

\begin{document}

\title{WarpagePINN: Thermal Warpage Prediction in Advanced Packaging via a Two-Stage Physics-Informed Neural Networks}

\author{Xinyu~Li,
        Min~Tang,~\IEEEmembership{Senior Member,~IEEE},
        Zeyu Sun,
        Wenxing~Zhu,
        Jianhua~Zhang,
        and Liang~Chen,~\IEEEmembership{Member,~IEEE} 
\thanks{This work was supported in part by National Natural Science Foundation of China under Grant 92473105 and 62504151; in part by National Key Research and Development Program of China under Grant 2025YFA1213000; and in part by State Key Laboratory of Radio Frequency Heterogeneous Integration (Open Scientific Research Program No. KF2024005). \it{(Corresponding author: Liang Chen.)}}
\thanks{X.~Li, J.~Zhang, and L.~Chen are with the School of
Microelectronics and Shanghai Key Laboratory of Chips and Systems for Intelligent Connected Vehicle, Shanghai University, Shanghai 201800, China (e-mail: lchenshu@shu.edu.cn).}
\thanks{M. Tang is with the State Key Laboratory of Radio Frequency Heterogeneous Integration, Shanghai Jiao Tong University, Shanghai, 200240, China.}
\thanks{Z. Sun is with the Institute of Microelectronics, Chinese Academy of Sciences, Beijing, 100029, China.}
\thanks{W. Zhu is with the Center for Discrete Mathematics and Theoretical
Computer Science, Fuzhou University, Fuzhou 350108, China.}
}


\maketitle
\begin{abstract}
  Thermal warpage has become a critical issue in advanced packaging, primarily caused by the mismatch in coefficients of thermal expansion (CTE) among heterogeneously integrated materials. However, only a limited number of studies have focused on developing computational methods for coupled thermal-warpage prediction in the chiplet. This paper proposes a two-stage physics-informed neural network (WarpagePINN) framework to compute both temperature profile and warpage deformation of chiplets. The neural networks are trained without relying on labeled datasets generated by conventional simulators. In the first stage, the temperature field is modeled using a Fourier series representation that inherently satisfies boundary conditions, and the network is trained solely through a loss function derived from the governing equation. In the second stage, a multilayer perceptron (MLP) is employed for warpage prediction, utilizing a novel hybrid supervisory strategy to optimize the energy-based loss function instead of residual loss. A parametric WarpagePINN is also developed to quantify uncertainties associated with the CTE. Numerical results show that the proposed WarpagePINN framework achieves excellent agreement with conventional finite element methods, with a mean absolute error (MAE) of 0.2 $\mu$m, while achieving a speedup of approximately 1000$\times$ in CTE parameterization studies.
\end{abstract}


\section{Introduction}
\label{sec:introduction}
\begin{figure}[ht]
\centering
\subfigure[]{
\includegraphics[width=0.81\linewidth]{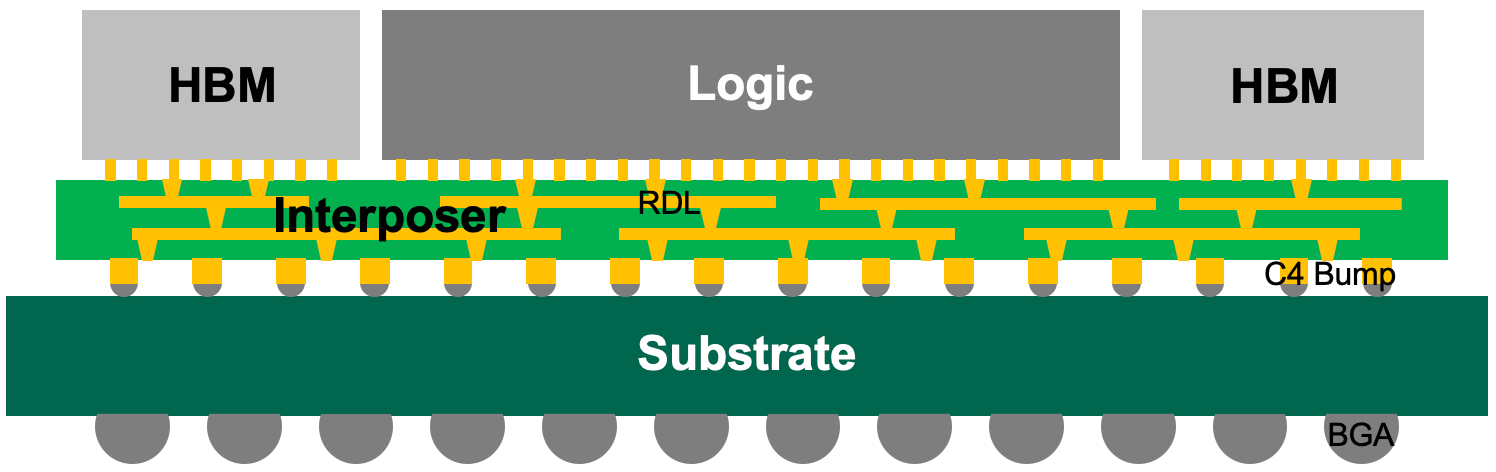}}
\subfigure[]{
\includegraphics[width=0.78\linewidth]{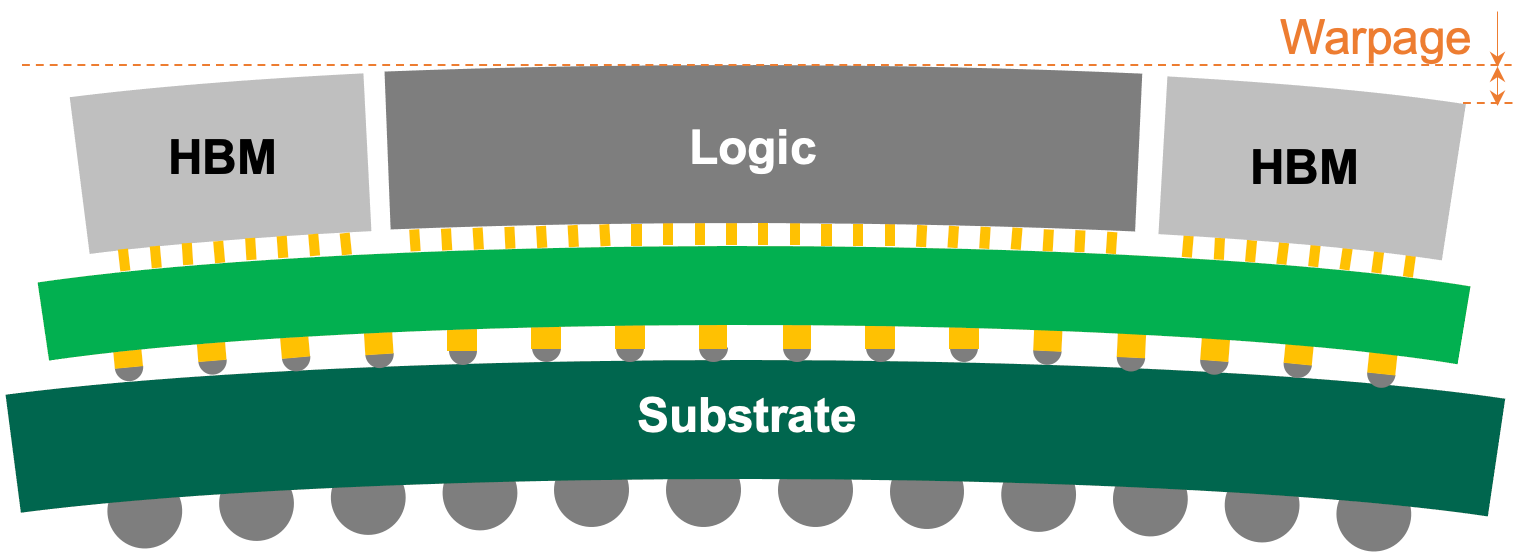}}
\caption{(a) Cross-sectional view of the TSMC CoWoS-R package. (b) Thermomechanical warpage of the package induced by heterogeneous integration and thermal expansion. }
\label{fig:1}
\end{figure}

\IEEEPARstart{A}{dvanced} packaging technology for multi-chiplet heterogeneous integration has become a promising alternative to transistor-level miniaturization for increasing integrated circuit density, effectively addressing the growing demands of high-performance computing (HPC) applications, including artificial intelligence (AI), large language models (LLMs) and autonomous vehicles~\cite{Su:ISSCC'23}. However, the use of disparate materials in such packages introduces a mismatch in their coefficients of thermal expansion (CTE)~\cite{Lau:TCPMT'17,Lau:TCPMT'18}. Moreover, continual increases in packaging dimensions, exemplified by state-of-the-art TSMC CoWoS technology with chip-on-wafer (CoW) sizes reaching 70 $\times$ 70 mm$^2$ and substrates up to 100 $\times$ 100 mm$^2$, have made thermal warpage a critical reliability issue~\cite{Hu:ECTC'25}, as shown in Fig.~\ref{fig:1}. Therefore, accurate prediction and mitigation of thermal warpage at early design stages are essential for modern heterogeneous integration systems~\cite{Chang:DT'25}.

Numerous modeling approaches have been established to predict warpage in electronic packages. Classical analytical models, such as laminate theory, Suhir's theory, and Timoshenko's beam model, offer valuable closed-form solutions for simplified structures and provide deep physical insights into the mechanics of stress and deformation~\cite{Timoshenko:JOSCA'25,Wang:TCPMT'21,Suhir:JRPC'93,Wang:EMAP'18,Tsai:TCPMT'20,Wei:TCPMT'12,Lo:ICCAD'24,Lo:TCAD'25}. For more complex, real-world geometries and material behaviors, numerical methods, primarily the Finite Element Method (FEM), have become the predominant tool~\cite{Zhu:DATE'25}. While FEM delivers high fidelity, its computational cost and meshing requirements can be prohibitive for rapid design iteration. This gap between the oversimplification of analytical models and the computational burden of FEM underscores the need for efficient and accurate alternatives.

Recently, the application of Artificial Intelligence (AI) for solving partial differential equations (PDEs) has witnessed significant breakthroughs, opening new paradigms for scientific computing~\cite{Lee:ECTC'25,Zhao:access'25,Lo:ICCAD'25,Raissi:JCP'19,Karniadakis:NRP'21,Chen:ICCAD'23}. These approaches can be broadly categorized into data-driven and physics-informed methods. While purely data-driven models rely on extensive simulation or experimental datasets to learn underlying patterns, they often struggle with generalization and physical consistency. In contrast, Physics-Informed Neural Networks (PINNs) have emerged as a powerful framework that seamlessly integrates the governing physical laws, expressed as PDEs, directly into the learning objective~\cite{Raissi:JCP'19,Karniadakis:NRP'21,Chen:ICCAD'23}. This allows PINNs to be trained with minimal labeled data while ensuring that the solutions adhere to the underlying physics, offering a promising path for accurate and efficient modeling of complex multi-physics problems.

In this paper, we propose a novel Two-Stage Physics-Informed Neural Network (WarpagePINN) framework to simultaneously predict the temperature distribution and warpage in chiplets by solving the governing heat conduction equation and solid mechanics equations.  Our contributions are outlined as follows:

\begin{itemlist}
\item We propose WarpagePINN, the first physics-informed neural network framework dedicated to fully coupled thermal–mechanical simulation of chiplet-based systems. The thermal field is obtained from the ARRR-PINN method~\cite{Zhou:DAC'25}, which employs a Fourier series representation and an adaptive sampling strategy to efficiently resolve temperature distributions.

\item The warpage field is characterized by a multilayer perceptron (MLP) trained within the PINN paradigm using two complementary loss formulations: an energy-based variational loss and a pure PDE-residual loss. A hybrid supervisory training strategy is introduced, where the energy-based loss drives parameter optimization while the residual loss guides the stopping criterion, thereby enhancing predictive accuracy and training stability.

\item A parameterized extension of WarpagePINN is developed to explicitly account for coefficients of thermal expansion (CTE). Numerical experiments demonstrate that the proposed approach achieves high predictive accuracy while delivering substantial computational speedup compared to conventional finite-element analysis.
\end{itemlist}




The remainder of this paper is organized as follows. Section~\ref{sec:RW} provides a review of related works. In Section~\ref{sec:elmsov}, we present the proposed WarpagePINN framework for the chiplet. Section~\ref{sec:experimentalresults} discusses the experimental results and analysis. Finally, conclusions are drawn in Section~\ref{sec:concl}.

\section{Relevant Work}
\label{sec:RW}
Based on the classical beam theory, Timoshenko pioneered the theoretical analysis of bi-metal thermostats~\cite{Timoshenko:JOSCA'25}. Building upon Timoshenko's model, Wang {\it et al.} developed an extended theoretical framework to characterize warpage evolution throughout various process stages. This model offers improved predictions for process-induced warpage and is particularly effective for multilayer film structures compared to the classical bimaterial approach~\cite{Wang:TCPMT'21}. For more complex material assemblies, Suhir derived a series of closed-form solutions for thermal stress and warpage in tri-material systems~\cite{Suhir:JRPC'93}. Tsai {\it et al.} later refined Suhir's solution by correcting the flexural rigidities of the plates to ensure consistency with the spherically bending hypothesis~\cite{Wang:EMAP'18}. Further extensions were made by Wang and Tsai, who presented solutions for assemblies with an extended base, incorporating the assumption of zero curvature across the entire base and accounting for the Poisson effect ~\cite{Tsai:TCPMT'20}. The classical laminated plate theory (CLPT) has also been widely applied to warpage analysis. For instance, Wei {\it et al.} employed CLPT to investigate the warpage behavior of printed wiring boards (PWBs) and PWB assemblies during the thermal convective reflow process. Closed-form solutions derived from CLPT were used to evaluate warpage at different locations across the samples~\cite{Wei:TCPMT'12}. More recently, Lo {\it et al.} proposed a highly efficient and accurate warpage model for advanced packages based on laminate theory, converting three-dimensional packaging structures into a two-dimensional thermo-solid coupling problem~\cite{Lo:ICCAD'24,Lo:TCAD'25}. However, this analytical approach still relies on a two-dimensional triangular finite element method to resolve warpage deformations in complex package geometries.

However, the aforementioned analytical studies are predominantly grounded in the central symmetry hypothesis and assume a fixed central point. These underlying assumptions significantly restrict the applicability of such models to modern asymmetric assemblies, which are increasingly common in 2.5D and 3D advanced packages integrating diverse components such as Memory, GPU, and High Bandwidth Memory (HBM). The analysis of asymmetric structures thus remains a key challenge for theoretical extension.
To address the limitations of analytical approaches, numerical methods based on Solid Mechanics theory have been widely adopted. Tools such as ANSYS and COMSOL utilize the Finite Element Method (FEM) to solve governing solid mechanics equations, providing accurate predictions of stress and warpage. To enhance computational efficiency, Zhu {\it et al.} proposed a model order reduction technique that exploits the periodicity of local fine structures, substantially reducing the number of degrees of freedom. When combined with sub-modeling techniques, this approach offers a flexible and efficient framework for accelerating finite element simulations~\cite{Zhu:DATE'25}. Despite these advancements, the high computational cost associated with full-scale FEM simulations, rooted in their physically rigorous formulation, often makes them impractical for direct application to multi-chiplet integration systems, where structural complexity and scale impose significant demands on computational resources.


Recently, machine learning methods for solving partial differential equations (PDEs) have demonstrated significant breakthroughs, offering novel data-driven approaches for addressing solid mechanics problems. In the context of package warpage prediction, several studies have explored such techniques. Lee {\it et al.} proposed a conditional generative adversarial network (cGAN)-based model to predict the global warpage distribution across the entire package surface~\cite{Lee:ECTC'25}. Similarly, Zhao {\it et al.} employed various artificial intelligence (AI) algorithms, including multilayer perceptron (MLP), extreme gradient boosting (XGB), and gradient boosting machine (GBM), to predict the warpage behavior of organic substrates under diverse material and structural configurations~\cite{Zhao:access'25}. Further advancing this direction, Lo {\it et al.} utilized a DeepONet-based operator learning framework to construct a model for predicting warpage in advanced packaging~\cite{Lo:ICCAD'25}. However, these purely data-driven deep learning methods rely heavily on labeled datasets generated from computationally expensive thermo-mechanical simulations, which demand substantial computational resources. To mitigate this dependency, physics-informed neural networks (PINNs) have been introduced. This approach integrates governing physical laws and boundary conditions directly into the loss function, thereby guiding the neural network training process without the need for extensive labeled data~\cite{Raissi:JCP'19,Karniadakis:NRP'21,Chen:ICCAD'23}. Despite its potential, no published work has yet applied the PINN framework to solve warpage problems in advanced packaging, indicating a promising direction for future research.


\section{The Proposed WarpagePINN Framework}
\label{sec:elmsov}

\begin{figure*}
    \centering
    \includegraphics[width=1\linewidth]{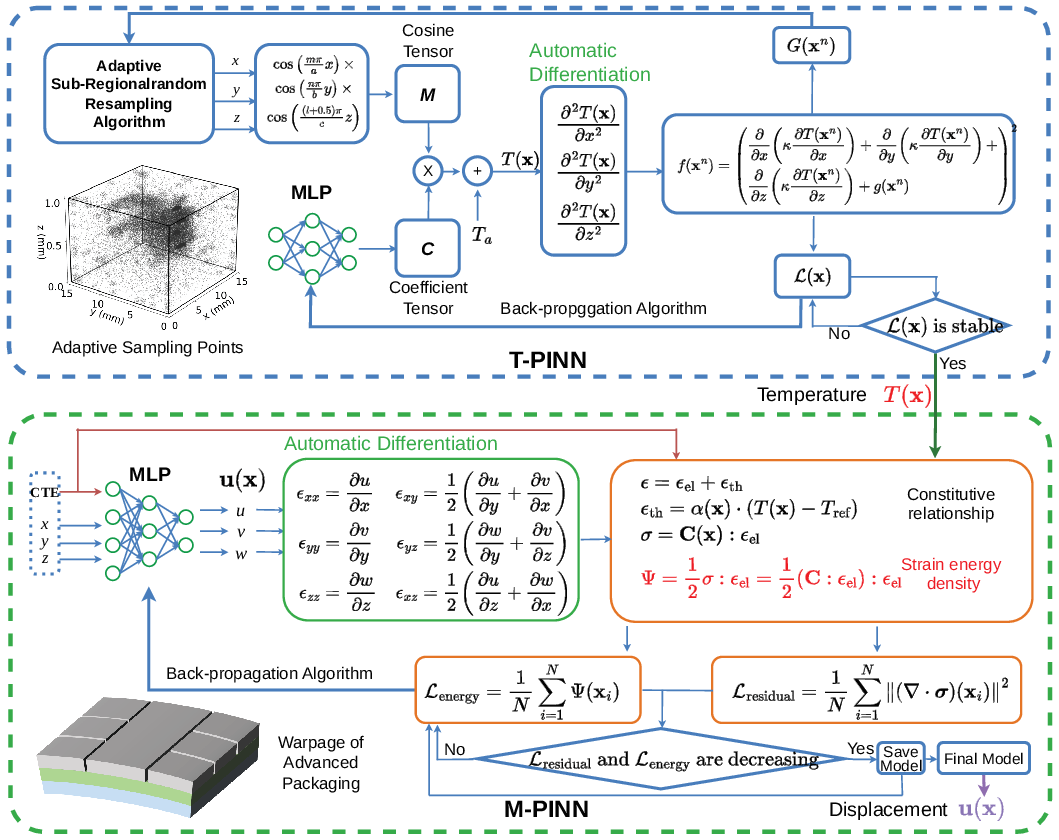}
    \caption{Framework of the proposed WarpagePINN for thermo-warpage prediction in advanced packaging.}
    \label{fig:placeholder}
\end{figure*}
The proposed WarpagePINN framework, illustrated in Fig.~\ref{fig:placeholder}, consists of two sequentially coupled PINNs: a Thermal PINN (T-PINN) for temperature field prediction and a Mechanical PINN (M-PINN) for warpage displacement computation. The T-PINN maps a given power density distribution to a steady-state temperature field. This temperature field then serves as a thermal load for the M-PINN, which predicts the resulting mechanical deformation. Both networks are trained without any labeled data from traditional solvers.


\subsection{Thermal Physics-Informed Neural Network}
The chip is modeled as a monolithic block with uniform thermal properties. Its bottom surface is maintained at a prescribed temperature \( T_a \), while all other surfaces are treated as adiabatic. Under these boundary conditions, the steady-state temperature distribution can be expressed by the following Fourier series expansion
\begin{equation}
\begin{split}
 T(\mathbf{x}) = \sum_{l=0}^{L} \sum_{m=0}^{M} \sum_{n=0}^{N} C_{lmn} \cos \left( \frac{m \pi}{a} x \right) \times \\ \cos \left( \frac{n \pi}{b} y \right) \cos \left( \frac{(l + 0.5) \pi}{c} z \right) + T_a,  
\end{split}
\label{eq:fsr}
\end{equation}
where \( a \), \( b \), and \( c \) denote the length, width, and height of the 3D-IC, respectively. The integers \( M \), \( N \), and \( L \) represent the truncation orders of the series along the \( x \), \( y \), and \( z \) directions, and \( C_{lmn} \) are the unknown coefficients corresponding to the basis function \( \cos \left( m \pi x/a  \right) \cos \left( n \pi y/b \right) \cos \left((l + 0.5) \pi z/c \right) \).

The choice of Fourier series as the output representation for T-PINN is motivated by several considerations. First, the basis functions $\cos(m\pi x/a)\cos(n\pi y/b)\cos((2l+1)\pi z/(2c))$ form a complete orthogonal set in $L^2(\Omega)$ and automatically satisfy the prescribed boundary conditions: adiabatic on side walls ($x=0,a$ and $y=0,b$) and bottom temperature fixed at $T_a$. This inductive bias significantly reduces the solution space that the network must explore. Second, Fourier bases are globally supported, enabling accurate approximation of smooth temperature fields with relatively few terms. The truncation orders $L,M,N$ are determined adaptively based on the smoothness of the power density distribution; for typical chiplet layouts with localized hot spots, we find that $L=M=N=20$ achieves a good balance between accuracy and computational efficiency. 

The PINN framework is then employed by incorporating the heat conduction equation into the loss function. Through minimization of this loss, the Fourier neural network (represented by Eq. \eqref{eq:fsr}) is trained to learn the coefficients \( C_{lmn} \). The loss function is defined as:
\begin{equation}
\mathcal{L}_{\text{loss}} = \left[ \kappa \left( \frac{\partial^2 T(\mathbf{x})}{\partial x^2} + \frac{\partial^2 T(\mathbf{x})}{\partial y^2} + \frac{\partial^2 T(\mathbf{x})}{\partial z^2} \right) + g(\mathbf{x}) \right]^2,
\end{equation}
where \( \kappa \) is the thermal conductivity and \( g(x,y,z) \) denotes the power density of the chiplet.

To enhance accuracy with a limited number of sampling points, the Adaptive Sub-regional Random Resampling (ASRR) method~\cite{Zhou:DAC'25} is adopted. This method adaptively adjusts the sampling distribution by increasing the density of points in regions with higher errors. In this work, ASRR-PINN is utilized to compute the temperature distribution for 3D-ICs.

\subsection{Mechanical Physics-Informed Neural Network}
\label{sec:MPINN}

\subsubsection{Linear Thermo-Elasticity Model}

In multi-chiplet heterogeneous integration, the deformation, though visually discernible, remains within the regime of small deformation theory (linear thermo-elasticity). This assumption linearizes the problem and offers computational efficiency. The displacement field $\mathbf{u}(\mathbf{x})$ and stress field $\boldsymbol{\sigma}$ satisfy three fundamental relations.

\noindent \textbf{Strain–Displacement Relationship.} Under the small deformation assumption, the Cauchy strain tensor $\boldsymbol{\epsilon}$ is linear in the displacement gradient:
\begin{equation}
\boldsymbol{\epsilon} = \frac{1}{2} \left( \nabla \mathbf{u} + (\nabla \mathbf{u})^T \right),
\end{equation}
where $\mathbf{u}(\mathbf{x}) = [u(\mathbf{x}), v(\mathbf{x}), w(\mathbf{x})]^T$ is the displacement vector. The components are
\[
\epsilon_{xx} = \frac{\partial u}{\partial x},\quad
\epsilon_{yy} = \frac{\partial v}{\partial y},\quad
\epsilon_{zz} = \frac{\partial w}{\partial z},
\]
\[\epsilon_{xy} = \epsilon_{yx} = \frac{1}{2}\left( \frac{\partial u}{\partial y} + \frac{\partial v}{\partial x} \right), \quad \dots\]

\noindent \textbf{Constitutive Law (Hooke’s Law).} The total strain $\boldsymbol{\epsilon}$ decomposes additively into an elastic part $\boldsymbol{\epsilon}_{\text{el}}$ and a thermal part $\boldsymbol{\epsilon}_{\text{th}}$:
\begin{equation}
\boldsymbol{\epsilon} = \boldsymbol{\epsilon}_{\text{el}} + \boldsymbol{\epsilon}_{\text{th}}.
\end{equation}
For an isotropic material, the thermal strain is
\begin{equation}
\boldsymbol{\epsilon}_{\text{th}} = \alpha(\mathbf{x}) \, (T(\mathbf{x}) - T_{\text{ref}}) \, \mathbf{I},
\end{equation}
where $\alpha(\mathbf{x})$ is the spatially varying coefficient of thermal expansion (CTE). The stress tensor $\boldsymbol{\sigma}$ is generated by the elastic strain via Hooke’s law:
\begin{equation}
\boldsymbol{\sigma} = \mathbf{C}(\mathbf{x}) : \boldsymbol{\epsilon}_{\text{el}},
\end{equation}
with $\mathbf{C}(\mathbf{x})$ the fourth‑order elasticity tensor, constructed from the spatially varying Young’s modulus $E(\mathbf{x})$ and Poisson’s ratio $\nu(\mathbf{x})$. The operator $:$ denotes the double‑dot product. 

For a general anisotropic material, $\mathbf{C}(\mathbf{x})$ has 21 independent components. However, in advanced packaging applications, materials are typically either isotropic (silicon, metals, polymers) or orthotropic (some composite substrates). For the common case of isotropic materials with spatially varying properties, the elasticity tensor simplifies to:

\begin{equation}
C_{ijkl}(\mathbf{x}) = \lambda(\mathbf{x}) \delta_{ij}\delta_{kl} + \mu(\mathbf{x}) (\delta_{ik}\delta_{jl} + \delta_{il}\delta_{jk})
\label{eq:isotropic_tensor}
\end{equation}
where $\lambda(\mathbf{x})$ and $\mu(\mathbf{x})$ are the spatially varying Lam\'e parameters, related to Young's modulus $E(\mathbf{x})$ and Poisson's ratio $\nu(\mathbf{x})$ by:

\begin{equation}
\lambda(\mathbf{x}) = \frac{E(\mathbf{x})\nu(\mathbf{x})}{(1+\nu(\mathbf{x}))(1-2\nu(\mathbf{x}))}, \quad
\mu(\mathbf{x}) = \frac{E(\mathbf{x})}{2(1+\nu(\mathbf{x}))}
\label{eq:lame_parameters}
\end{equation}
In Voigt notation (which maps symmetric tensor pairs to vector indices: $11 \rightarrow 1$, $22 \rightarrow 2$, $33 \rightarrow 3$, $23 \rightarrow 4$, $13 \rightarrow 5$, $12 \rightarrow 6$), the constitutive relation for an isotropic material can be written in matrix form as:

\begin{figure*}[htb]
\begin{align}
\begin{Bmatrix}
\sigma_{xx} \\ \sigma_{yy} \\ \sigma_{zz} \\ \sigma_{yz} \\ \sigma_{xz} \\ \sigma_{xy}
\end{Bmatrix} = 
\frac{E(\mathbf{x})}{(1+\nu(\mathbf{x}))(1-2\nu(\mathbf{x}))}
\begin{bmatrix}
1-\nu(\mathbf{x}) & \nu(\mathbf{x}) & \nu(\mathbf{x}) & 0 & 0 & 0 \\
\nu(\mathbf{x}) & 1-\nu(\mathbf{x}) & \nu(\mathbf{x}) & 0 & 0 & 0 \\
\nu(\mathbf{x}) & \nu(\mathbf{x}) & 1-\nu(\mathbf{x}) & 0 & 0 & 0 \\
0 & 0 & 0 & \frac{1-2\nu(\mathbf{x})}{2} & 0 & 0 \\
0 & 0 & 0 & 0 & \frac{1-2\nu(\mathbf{x})}{2} & 0 \\
0 & 0 & 0 & 0 & 0 & \frac{1-2\nu(\mathbf{x})}{2}
\end{bmatrix}
\begin{Bmatrix}
(\epsilon_{\text{el}})_{xx} \\ (\epsilon_{\text{el}})_{yy} \\ (\epsilon_{\text{el}})_{zz} \\ 2(\epsilon_{\text{el}})_{yz} \\ 2(\epsilon_{\text{el}})_{xz} \\ 2(\epsilon_{\text{el}})_{xy}
\end{Bmatrix}
\label{eq:constitutive_matrix}
\end{align}
\end{figure*}

For heterogeneous structures with multiple material layers, $E(\mathbf{x})$, $\nu(\mathbf{x})$, and $\alpha(\mathbf{x})$ are piecewise constant functions, with sharp transitions at layer interfaces. In practice, we use a smooth approximation (e.g., hyperbolic tangent functions) to avoid numerical instabilities during automatic differentiation.

\noindent \textbf{Force Equilibrium.} In the static case with negligible body forces, equilibrium requires
\begin{equation}
\nabla \cdot \boldsymbol{\sigma} = \mathbf{0}.
\label{eq:strong_form}
\end{equation}
In component form, this represents three partial differential equations:

\begin{equation}
\frac{\partial \sigma_{xx}}{\partial x} + \frac{\partial \sigma_{xy}}{\partial y} + \frac{\partial \sigma_{xz}}{\partial z} = 0
\end{equation}
\begin{equation}
\frac{\partial \sigma_{xy}}{\partial x} + \frac{\partial \sigma_{yy}}{\partial y} + \frac{\partial \sigma_{yz}}{\partial z} = 0
\end{equation}
\begin{equation}
\frac{\partial \sigma_{xz}}{\partial x} + \frac{\partial \sigma_{yz}}{\partial y} + \frac{\partial \sigma_{zz}}{\partial z} = 0
\end{equation}

On free boundaries (all external surfaces of the chiplet), the traction-free condition applies:

\begin{equation}
\mathbf{n} \cdot \boldsymbol{\sigma} = \mathbf{0} \quad \text{on} \quad \partial \Omega
\label{eq:traction_free}
\end{equation}
where $\mathbf{n}$ is the outward unit normal vector to the boundary $\partial \Omega$.

\subsubsection{Parameterized PINN Framework with a Hybrid Supervisory Training Strategy}

A multilayer perceptron (MLP) is employed to represent the displacement field. To accommodate heterogeneous material properties, the MLP takes both the spatial coordinates and the local CTE value as inputs:
\begin{equation}
\mathbf{u}(\mathbf{x}) = \text{MLP}\big(\mathbf{x}, \alpha(\mathbf{x})\big).
\end{equation}
This parameterization allows the network to learn the mapping from location and material property to displacement. If desired, other material fields (e.g., $E(\mathbf{x})$, $\nu(\mathbf{x})$) can be added as inputs or themselves represented by auxiliary networks. Leveraging automatic differentiation, all required spatial derivatives are computed, enabling the formulation of loss functions that enforce the thermo‑elasticity laws.

\textbf{Weak-Form Energy-Based Loss.}To improve training stability, we exploit the weak form of the linear thermo‑elasticity model. According to the principle of minimum potential energy, the stable equilibrium displacement field minimizes the total potential energy $\Pi$. In the absence of external forces, $\Pi$ reduces to the total strain energy:
\begin{equation}
\Pi = \int_{\Omega} \Psi(\boldsymbol{\epsilon}_{\text{el}}) \, dV,
\end{equation}
where $\Psi$ is the strain energy density. For a linear elastic material,
\begin{equation}
\Psi = \frac{1}{2} \boldsymbol{\sigma} : \boldsymbol{\epsilon}_{\text{el}} = \frac{1}{2} (\mathbf{C} : \boldsymbol{\epsilon}_{\text{el}}) : \boldsymbol{\epsilon}_{\text{el}}.
\label{eq:energy_density}
\end{equation}
In terms of engineering constants for isotropic materials, this expands to:

\begin{equation}
\Psi = \frac{1}{2} \left[ \lambda (\text{tr}(\boldsymbol{\epsilon}_{\text{el}}))^2 + 2\mu \, \text{tr}(\boldsymbol{\epsilon}_{\text{el}}^2) \right]
\label{eq:energy_density_isotropic}
\end{equation}
or, in component form:

\begin{equation}
\begin{split}
&\Psi = \frac{1}{2} [ \lambda (\epsilon_{xx} + \epsilon_{yy} + \epsilon_{zz})^2 +\\& 2\mu (\epsilon_{xx}^2 + \epsilon_{yy}^2 + \epsilon_{zz}^2 + 2\epsilon_{xy}^2 + 2\epsilon_{xz}^2 + 2\epsilon_{yz}^2) ]
\end{split}
\label{eq:energy_density_components}
\end{equation}
The system is fully determined when $\Pi = 0$, leading to the energy‑based loss function:
\begin{equation}
\mathcal{L}_{\text{energy}} = \frac{1}{N} \sum_{i=1}^{N} \Psi\big(\mathbf{x}_i, \alpha(\mathbf{x}_i)\big),
\label{eq:loss_energy}
\end{equation}
with $N$ collocation points in the domain $\Omega$. Minimizing $\mathcal{L}_{\text{energy}}$ is equivalent to solving the variational problem; it naturally enforces Neumann boundary conditions and acts as a smoothing operator, making optimization with gradient‑based methods (e.g., Adam, LBFGS) more robust.

\textbf{Strong-Form Residual Loss.} Alternatively, one may directly penalize the residuals of the strong‑form equations:
\begin{equation}
\begin{split}
\mathcal{L}_{\text{residual}} = \frac{1}{N} \sum_{i=1}^{N} \big\| (\nabla \cdot \boldsymbol{\sigma})(\mathbf{x}_i, \alpha(\mathbf{x}_i)) \big\|^2 \\+ \frac{1}{N_b} \sum_{i=1}^{N_b} \big\| \mathbf{n} \cdot \boldsymbol{\sigma}(\mathbf{x}_i) \big\|^2.
\end{split}
\label{eq:loss_residual}
\end{equation}
However, $\mathcal{L}_{\text{residual}}$ alone often leads to difficult optimization and poor convergence.

\noindent \textbf{Hybrid Supervisory Training Strategy.}
To combine the advantages of both loss formulations, we introduce a hybrid supervisory strategy that monitors both $\mathcal{L}_{\text{energy}}$ and $\mathcal{L}_{\text{residual}}$ during training. The network is trained by minimizing $\mathcal{L}_{\text{energy}}$, but the model is saved only when both losses achieve new minimum values. This ensures that reductions in system energy are accompanied by improved fidelity to the governing partial differential equations. The procedure is summarized in Algorithm \ref{alg:hybrid}.

\begin{algorithm}[ht]
\caption{Hybrid Supervisory Training Strategy}
\label{alg:hybrid}
\KwIn{Training points $\{\mathbf{x}_i\}_{i=1}^N$, CTE values $\{\alpha(\mathbf{x}_i)\}_{i=1}^N$, initial network parameters $\boldsymbol{\theta}$}
\KwOut{Optimized model parameters}
$\mathcal{L}_{\text{energy}}^{\text{best}} \leftarrow \infty$, $\mathcal{L}_{\text{residual}}^{\text{best}} \leftarrow \infty$\;
\For{epoch $k = 1$ \KwTo $K$}{
    Compute $\mathcal{L}_{\text{energy}}(\boldsymbol{\theta})$ and $\mathcal{L}_{\text{residual}}(\boldsymbol{\theta})$ on all points\;
    Update $\boldsymbol{\theta}$ by minimizing $\mathcal{L}_{\text{energy}}$ (e.g., using Adam or LBFGS)\;
    \If{$\mathcal{L}_{\text{energy}}(\boldsymbol{\theta}) < \mathcal{L}_{\text{energy}}^{\text{best}}$ \textbf{and} $\mathcal{L}_{\text{residual}}(\boldsymbol{\theta}) < \mathcal{L}_{\text{residual}}^{\text{best}}$}{
        $\mathcal{L}_{\text{energy}}^{\text{best}} \leftarrow \mathcal{L}_{\text{energy}}(\boldsymbol{\theta})$\;
        $\mathcal{L}_{\text{residual}}^{\text{best}} \leftarrow \mathcal{L}_{\text{residual}}(\boldsymbol{\theta})$\;
        Save current model parameters $\boldsymbol{\theta}$\;
    }
}
\end{algorithm}

This hybrid approach yields a model that simultaneously satisfies the energy minimum and the equilibrium equations, providing accurate displacement predictions while maintaining physical consistency.

\section{Experimental Results and Discussions}
\label{sec:experimentalresults}
\begin{figure*}[htb]
\centering
\subfigure[]{
\includegraphics[width=0.34\linewidth]{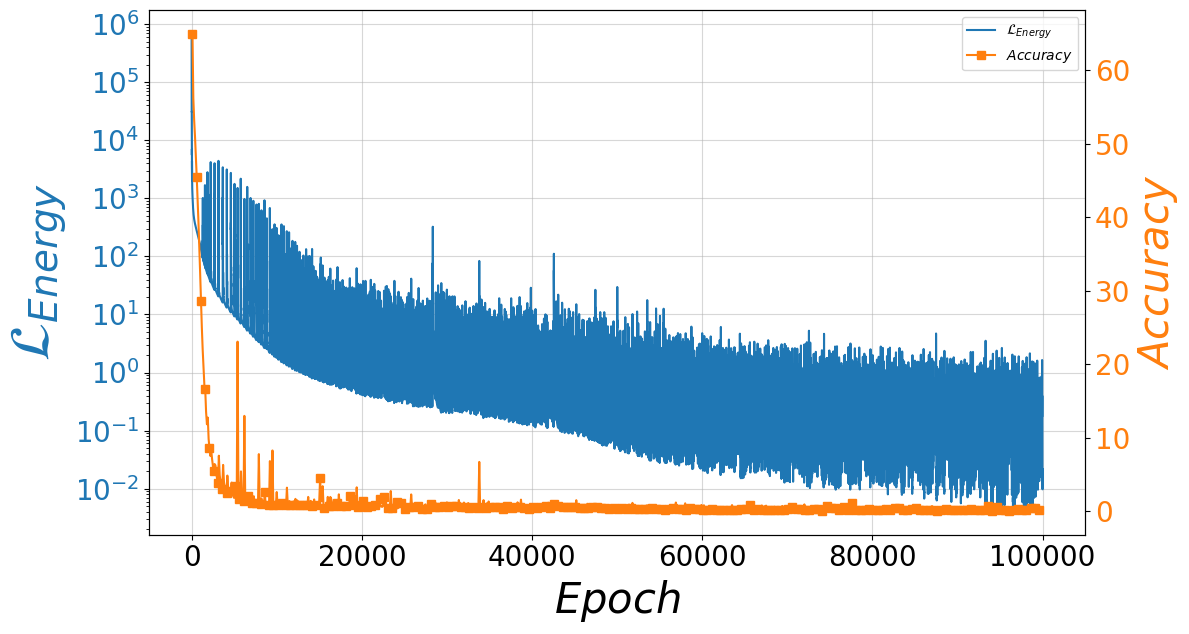}\label{fig:le}}
\subfigure[]{
\includegraphics[width=0.34\linewidth]{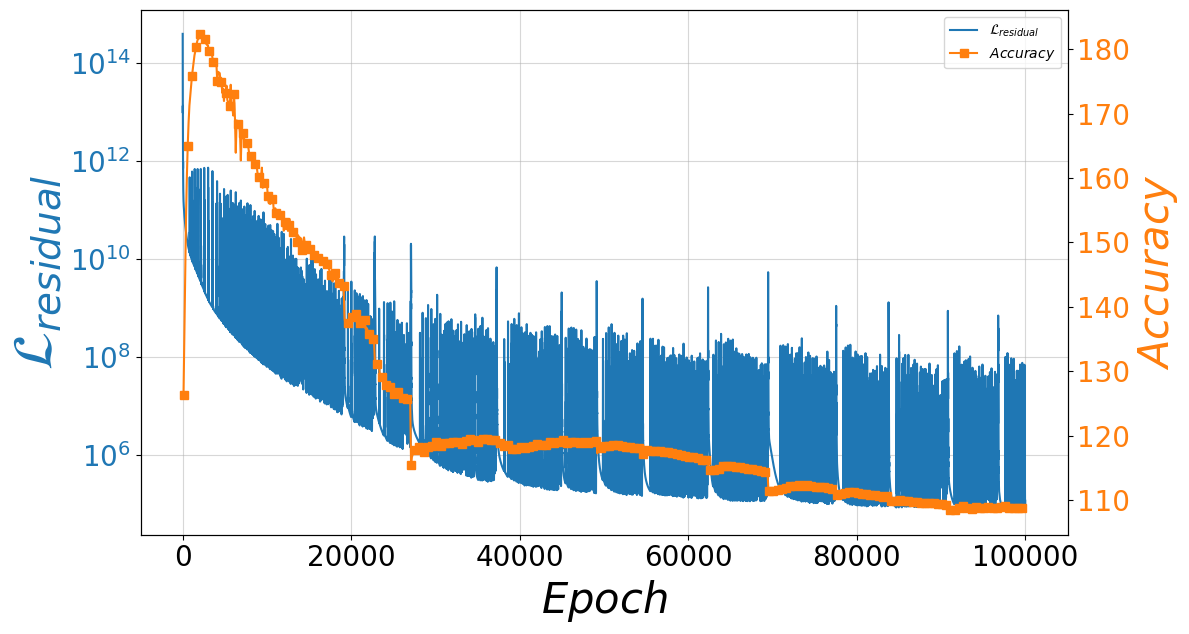}\label{fig:ls}}
\subfigure[]{
\includegraphics[width=0.28\linewidth]{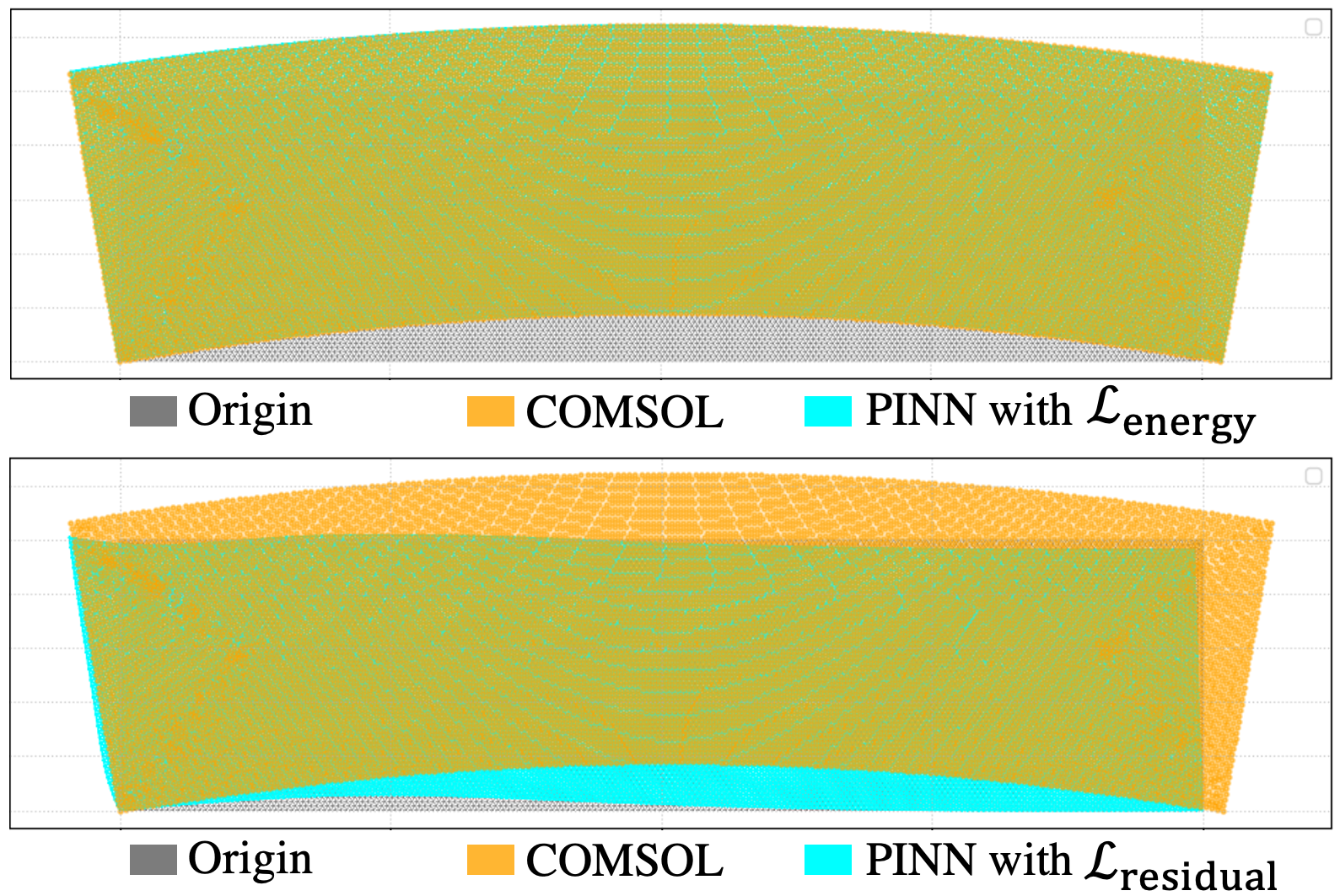}\label{fig:lesc}}
\caption{(a) Training history of the energy‑based loss $\mathcal{L}_{\text{energy}}$.
(b) Training history of the PDE‑based loss $\mathcal{L}_{\text{residual}}$.
(c) Comparison of results obtained with the two loss functions.}
\label{fig:lossfunction}
\end{figure*}
In this section, heterogeneous chip structures with diverse layouts are adopted to evaluate the accuracy and efficiency of the proposed WarpagePINN framework. First, the superiority of the energy-based loss function compared to the residual-based formulation is validated. Second, the detailed workflow of WarpagePINN is demonstrated through a TSMC CoWoS package case study. Third, the predictive performance of both T-PINN and M-PINN is evaluated across eight layout cases. Finally, a parametric study on the CTE is conducted using WarpagePINN.

The proposed WarpagePINN framework is implemented using the PyTorch platform. Ground truth solutions for the temperature and displacement fields are obtained using COMSOL Multiphysics. All simulations are performed on a Linux server equipped with an Intel Xeon W9-3495X 1.9 GHz CPU, two NVIDIA RTX A6000 48 GB GPUs, and two NVIDIA RTX 4090 24 GB GPUs.

\subsection{Comparison of $\mathcal{L}_{\text{energy}}$ and $\mathcal{L}_{\text{residual}}$}

To demonstrate the advantages of the proposed energy-based loss function $\mathcal{L}_{\text{energy}}$, which is derived from the variational weak form of the governing partial differential equations (PDEs), a comparative analysis was conducted against the conventional residual-based loss $\mathcal{L}_{\text{residual}}$ based on the strong-form formulation. Both methodologies were evaluated using a benchmark 2D chip warpage problem.

As illustrated in Fig. \ref{fig:le}, the energy-based loss exhibits steady convergence, ultimately achieving high predictive fidelity. In contrast, while the residual-based loss also reaches a stable state, it yields significantly higher prediction errors, as shown in Fig. \ref{fig:ls}. A spatial visualization in Fig. \ref{fig:lesc} further confirms that the displacement field predicted via $\mathcal{L}_{\text{energy}}$ aligns closely with the COMSOL reference solution, whereas the $\mathcal{L}_{\text{residual}}$ approach fails to capture the correct deformation profile accurately.

Furthermore, the computational efficiency of the energy-based framework is markedly superior; the training process required only 230 s, representing a substantial reduction compared to the 604 s required by the residual-based approach. These results validate that the energy-based loss offers both enhanced numerical accuracy and significantly lower computational overhead. Consequently, the energy-based formulation is adopted for all subsequent investigations in this study.

The superior performance of the energy-based loss $\mathcal{L}_{\text{energy}}$ over the residual-based loss $\mathcal{L}_{\text{residual}}$ can be understood through the lens of optimization theory. Let $\mathcal{F}$ denote the space of admissible displacement fields. The residual loss involves second-order derivatives of the network output, leading to a highly non-convex and ill-conditioned landscape:
\begin{equation}
\mathcal{L}_{\text{residual}}(\mathbf{u}) = \|\nabla \cdot \boldsymbol{\sigma}(\mathbf{u})\|^2 = \|\mathcal{P}(\mathbf{u})\|^2,
\end{equation}
where $\mathcal{P}$ is a fourth-order differential operator. The Hessian of $\mathcal{L}_{\text{residual}}$ contains high-frequency components that can cause gradient vanishing/exploding during backpropagation.

In contrast, the energy-based loss involves only first-order derivatives:
\begin{equation}
\mathcal{L}_{\text{energy}}(\mathbf{u}) = \int_\Omega \Psi(\boldsymbol{\epsilon}_{\text{el}}) dV = \frac{1}{2} \int_\Omega \boldsymbol{\epsilon}_{\text{el}} : \mathbf{C} : \boldsymbol{\epsilon}_{\text{el}} dV.
\end{equation}
This quadratic form in the strain (first derivatives of $\mathbf{u}$) yields a smoother landscape with better-conditioned Hessian. 
This analysis explains the 2.6$\times$ training speedup and superior accuracy observed with $\mathcal{L}_{\text{energy}}$.

\subsection{Thermal warpage predictions of a 3D CoWoS package}

\begin{figure}[ht]
\centering
\subfigure[]{
\includegraphics[width=0.77\linewidth]{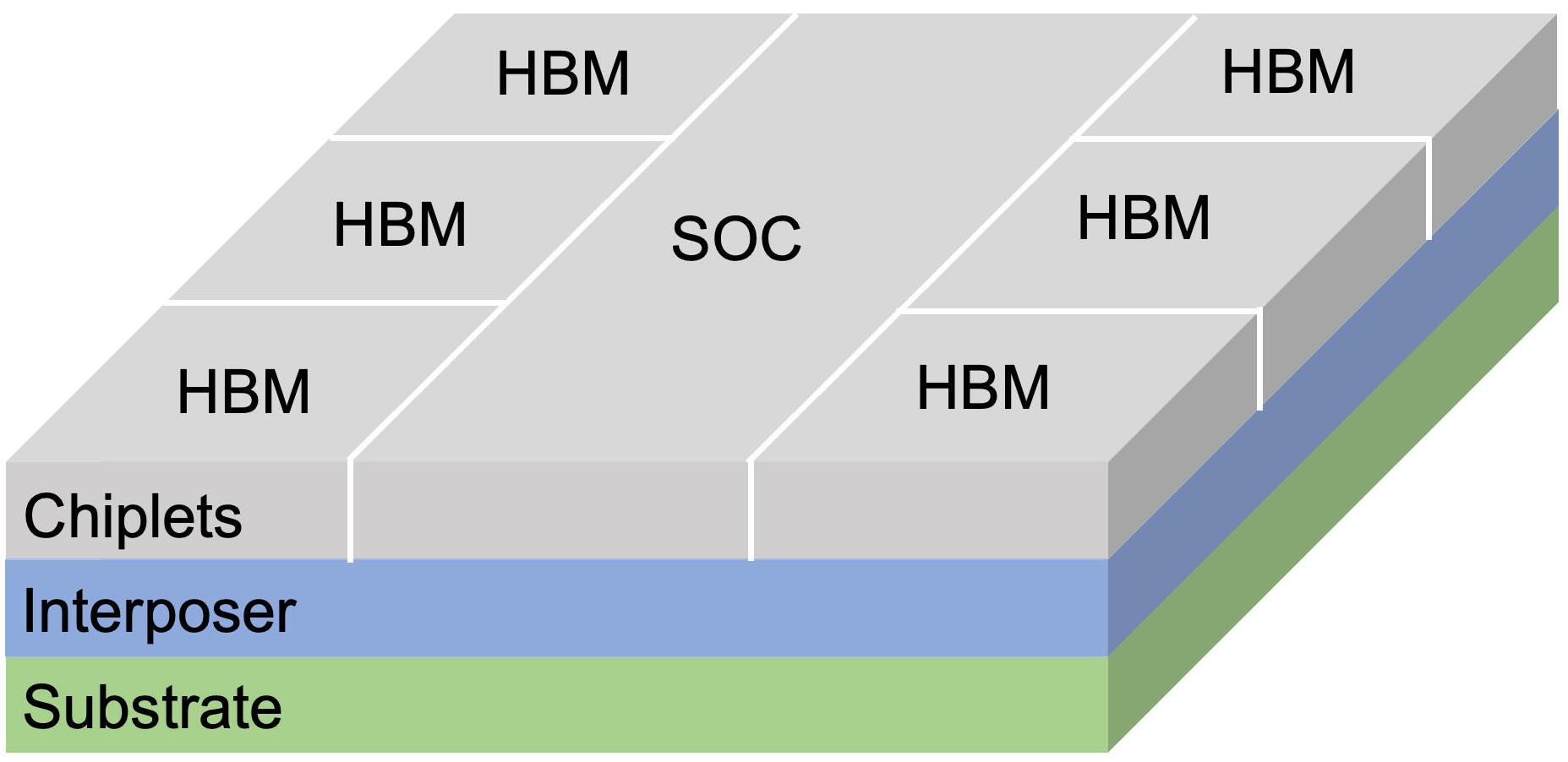}\label{fig:3ds}}
\subfigure[]{
\includegraphics[width=0.47\linewidth]{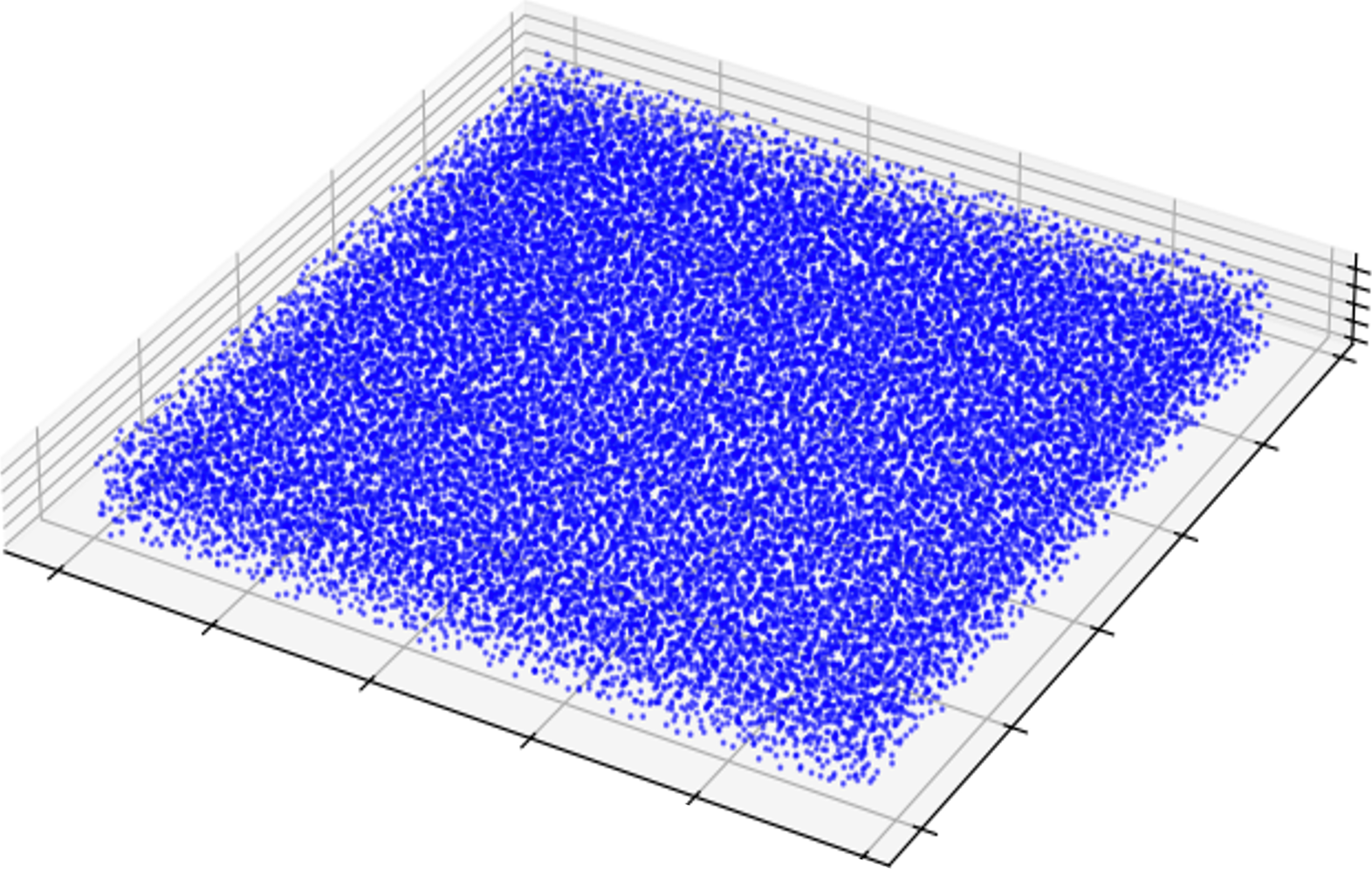}\label{fig:mesh}}
\subfigure[]{
\includegraphics[width=0.47\linewidth]{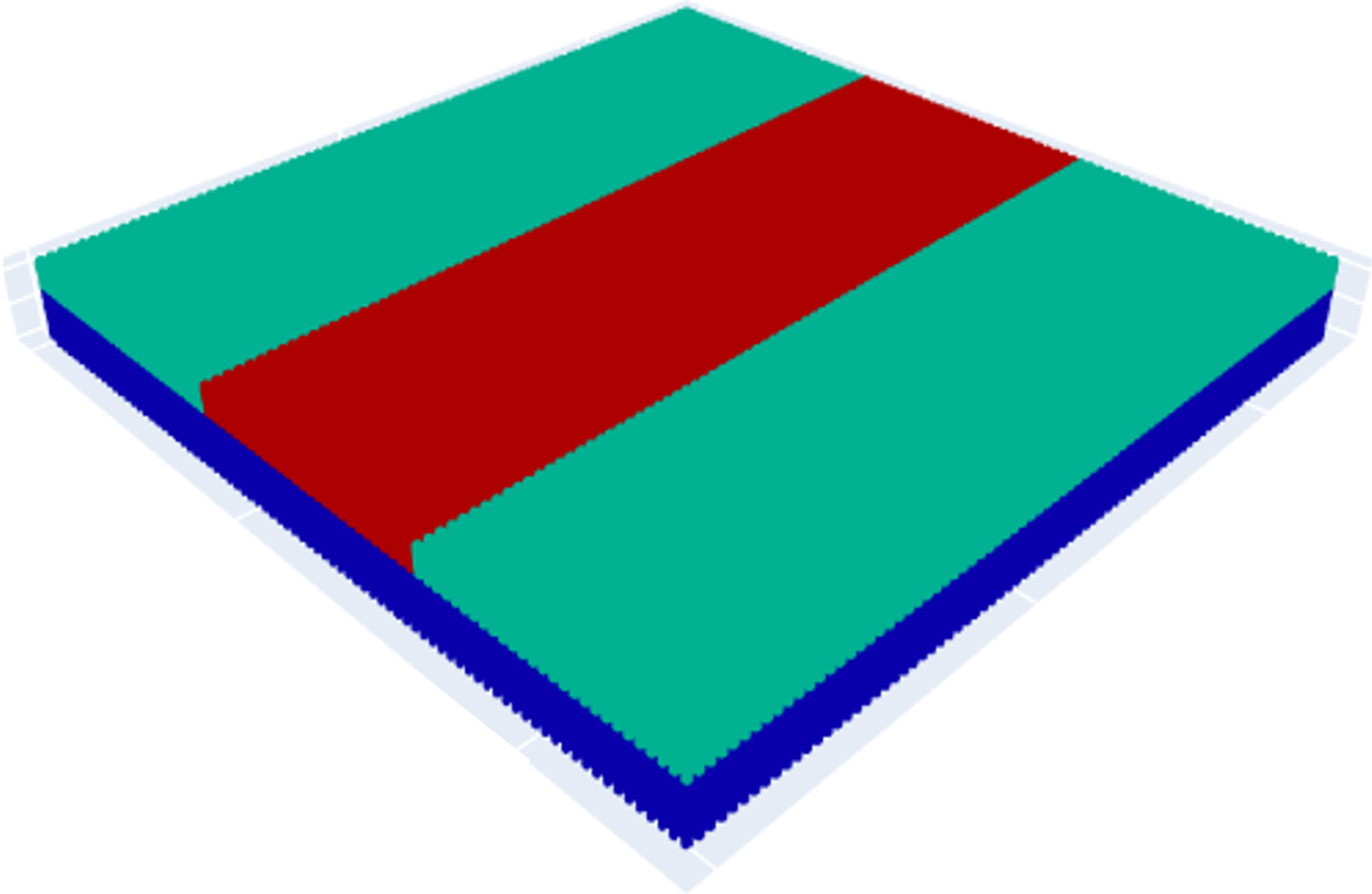}\label{fig:power}}
\subfigure[]{
\includegraphics[width=0.47\linewidth]{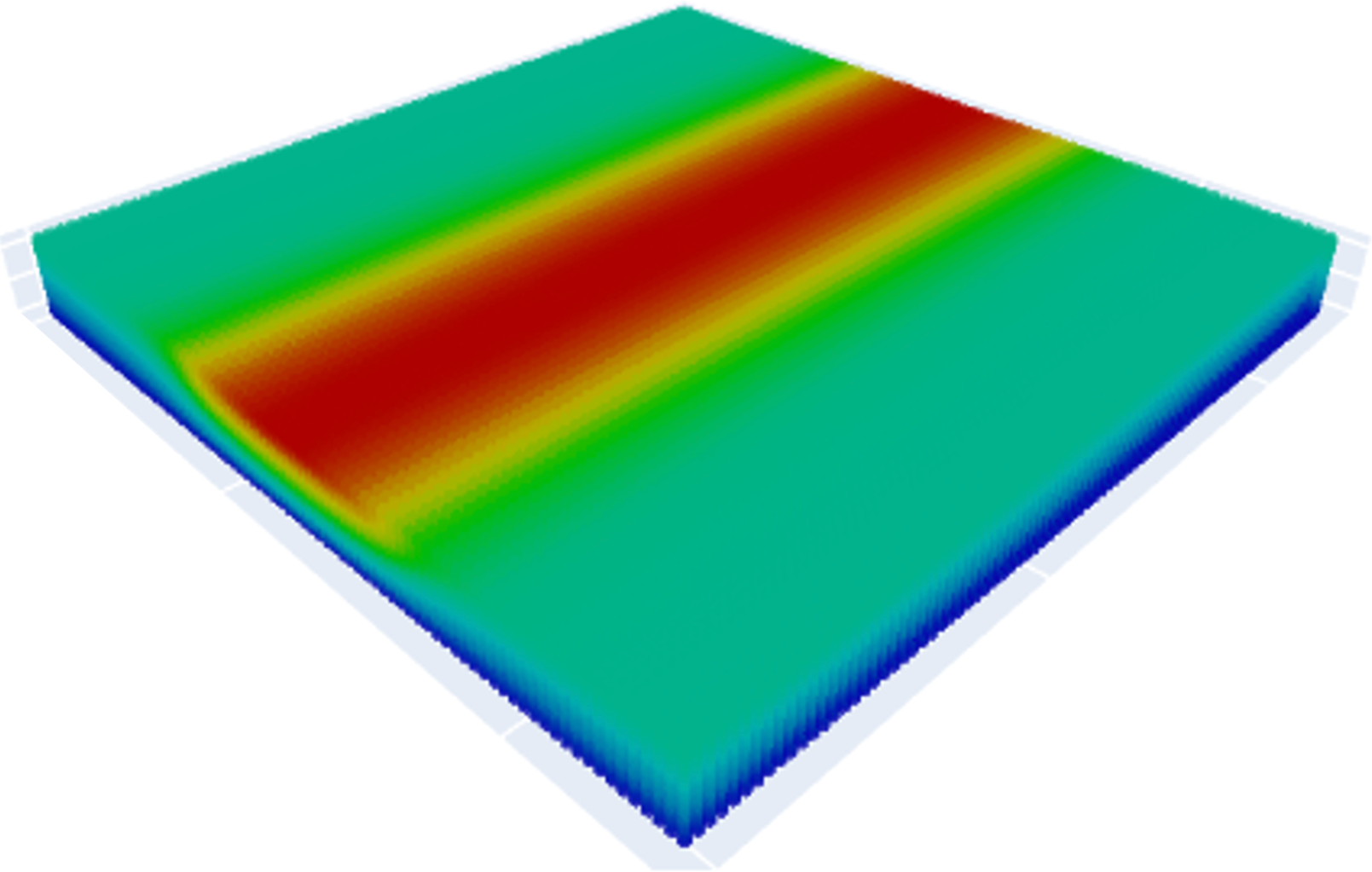}\label{fig:temp}}
\subfigure[]{
\includegraphics[width=0.47\linewidth]{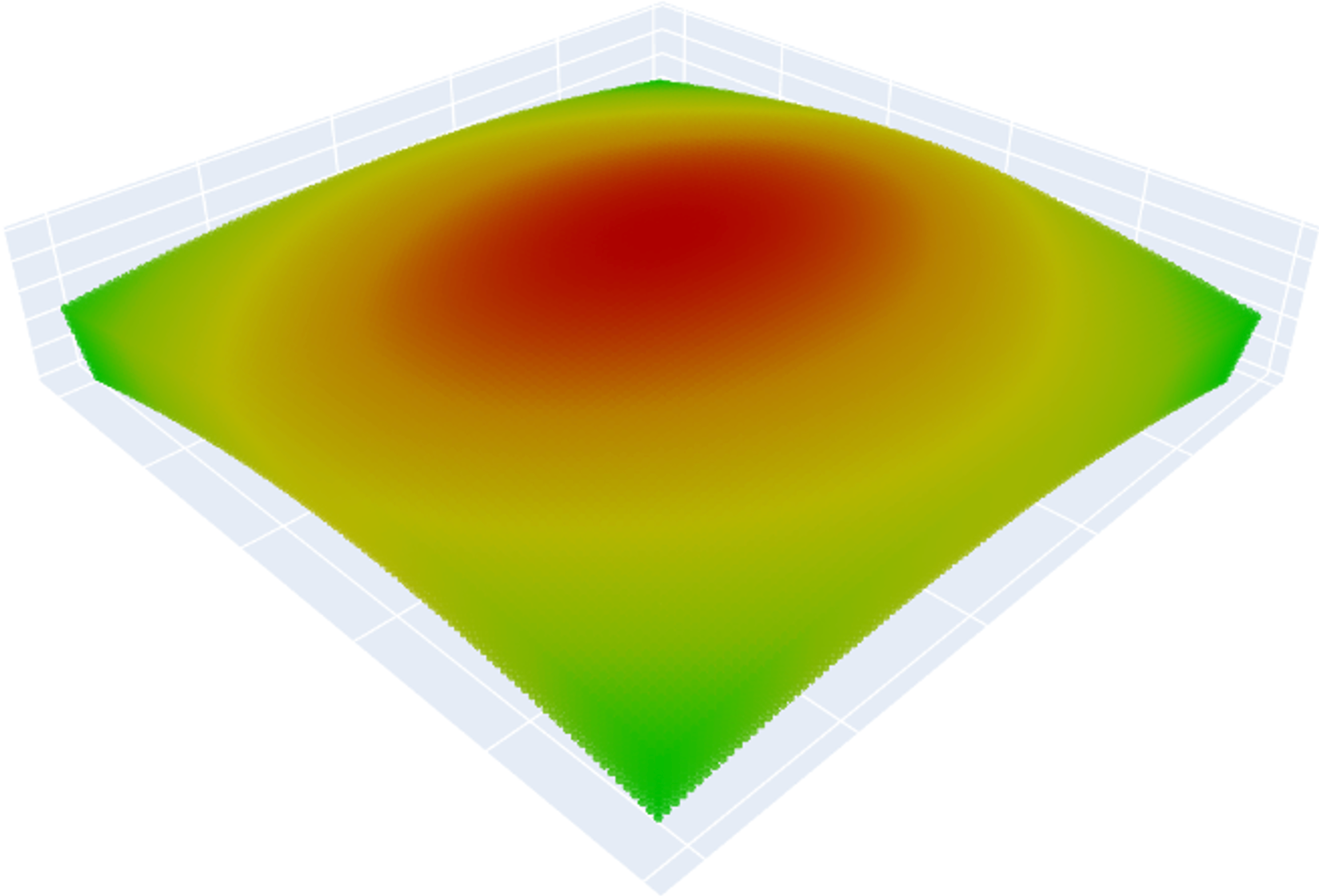}\label{fig:disp}}
\subfigure[]{
\includegraphics[width=0.47\linewidth]{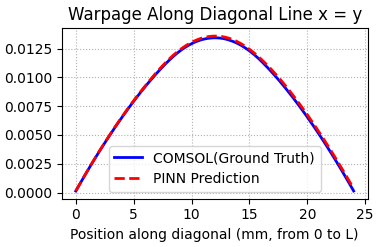}\label{fig:line1}}
\subfigure[]{
\includegraphics[width=0.47\linewidth]{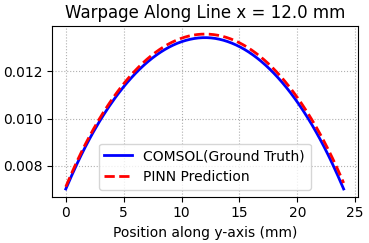}\label{fig:line2}}
\caption{(a) Schematic of a simplified TSMC CoWoS-R advanced package. (b) Sampling points. (c) Power density. (d) Temperature profile. (e) Warpage (magnified by a factor of 100). Comparison of two methods (f) along the diagonal line ($x = y$) and (g) along the line $x = 12\text{ mm}$.}
\label{fig:cowos}
\end{figure}


\begin{figure*}[htbp]
\centering

\setlength{\tabcolsep}{0pt}        
\renewcommand{\arraystretch}{0}    
\setlength{\extrarowheight}{0pt}  

\begin{tabular}{@{}cccccccc@{}}

\raisebox{0.9\imgheighth}{\textbf{Case 1}} &
\includegraphics[height=\imgheight]{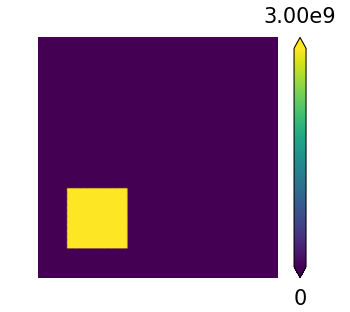} &
\hspace{-4pt}\includegraphics[height=\imgheight]{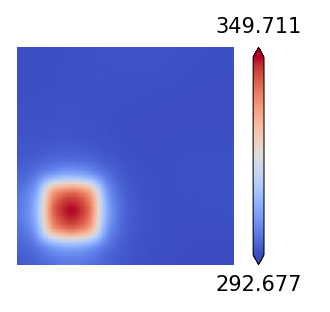} &
\hspace{-2pt}\includegraphics[height=\imgheight]{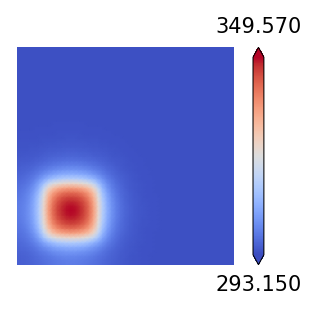} &
\hspace{-2pt}\includegraphics[height=\imgheight]{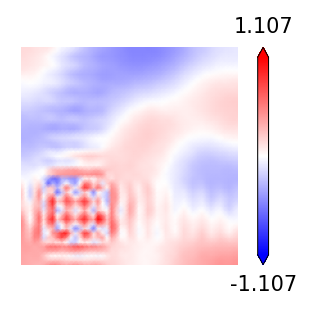} &
\hspace{-2pt}\includegraphics[height=\imgheight]{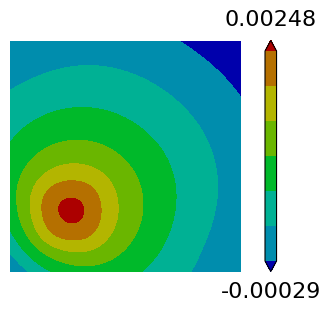} &
\hspace{-2pt}\includegraphics[height=\imgheight]{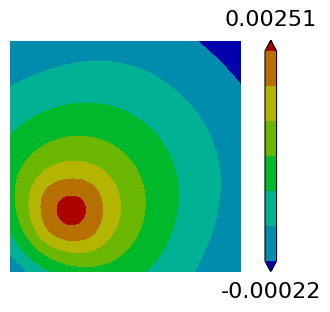} &
\hspace{-2pt}\includegraphics[height=\imgheight]{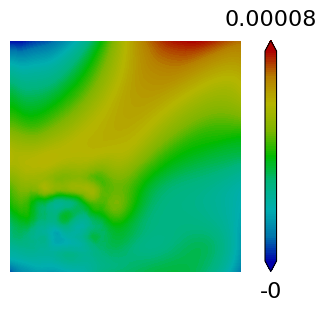} \\[-1pt]

\raisebox{0.9\imgheighth}{\textbf{Case 2}} &
\includegraphics[height=\imgheight]{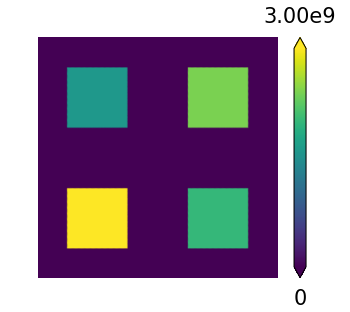} &
\hspace{-2pt}\includegraphics[height=\imgheight]{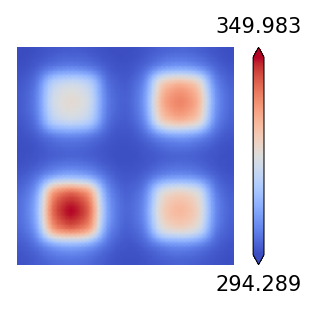} &
\hspace{-2pt}\includegraphics[height=\imgheight]{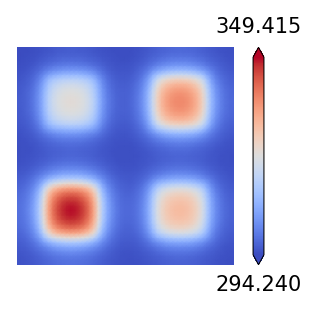} &
\hspace{-2pt}\includegraphics[height=\imgheight]{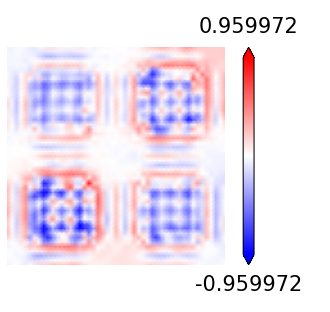} &
\hspace{-2pt}\includegraphics[height=\imgheight]{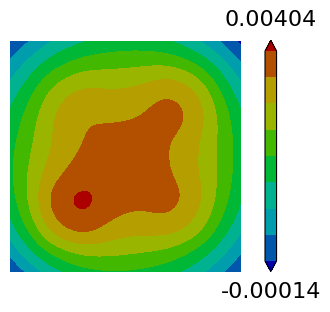} &
\hspace{-2pt}\includegraphics[height=\imgheight]{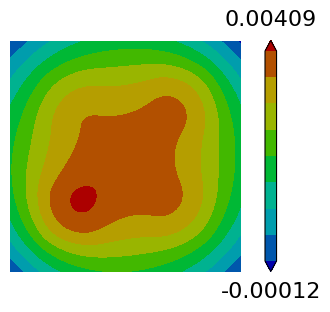} &
\hspace{-2pt}\includegraphics[height=\imgheight]{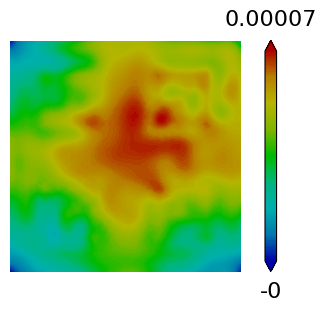} \\[-1pt]

\raisebox{0.9\imgheighth}{\textbf{Case 3}} &
\includegraphics[height=\imgheight]{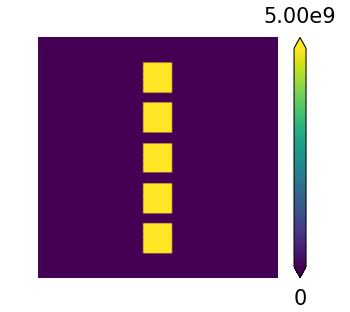} &
\hspace{-2pt}\includegraphics[height=\imgheight]{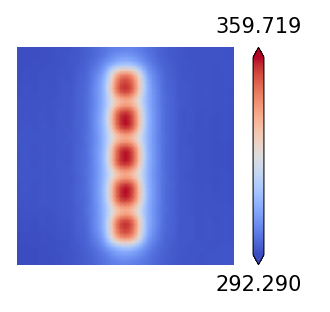} &
\hspace{-2pt}\includegraphics[height=\imgheight]{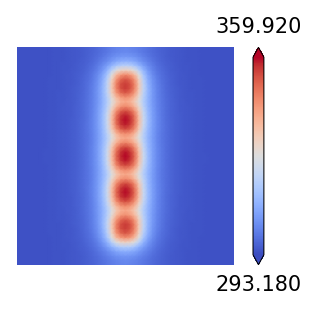} &
\hspace{-2pt}\includegraphics[height=\imgheight]{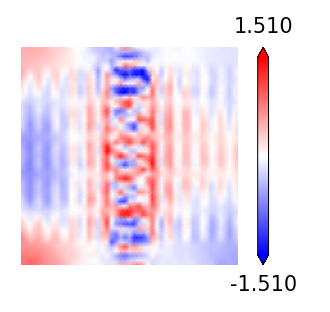} &
\hspace{-2pt}\includegraphics[height=\imgheight]{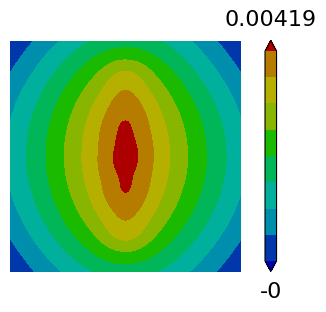} &
\hspace{-2pt}\includegraphics[height=\imgheight]{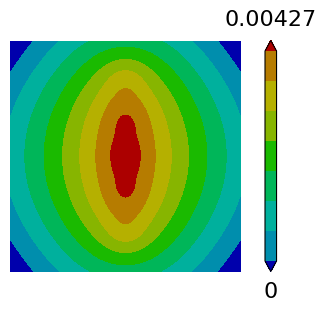} &
\hspace{-2pt}\includegraphics[height=\imgheight]{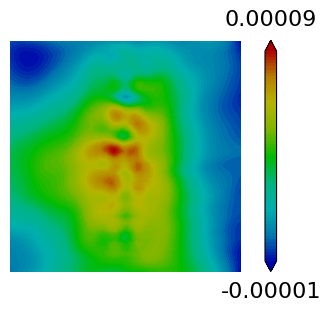} \\[-1pt]

\raisebox{0.9\imgheighth}{\textbf{Case 4}} &
\includegraphics[height=\imgheight]{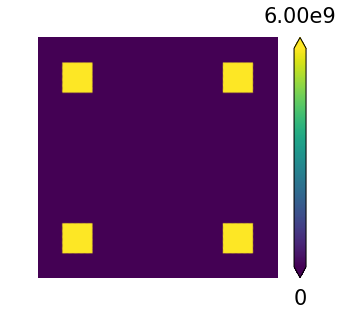} &
\hspace{-2pt}\includegraphics[height=\imgheight]{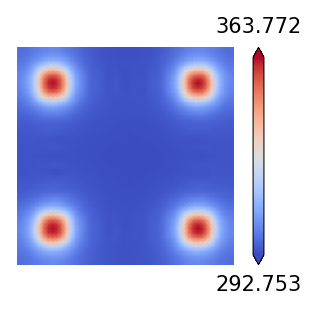} &
\hspace{-2pt}\includegraphics[height=\imgheight]{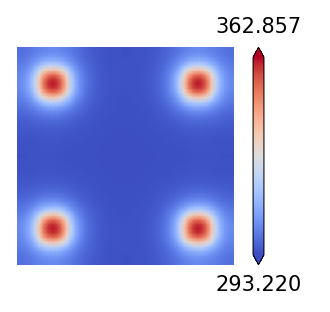} &
\hspace{-2pt}\includegraphics[height=\imgheight]{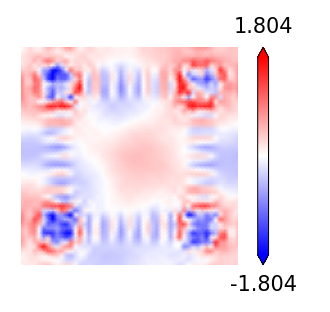} &
\hspace{-2pt}\includegraphics[height=\imgheight]{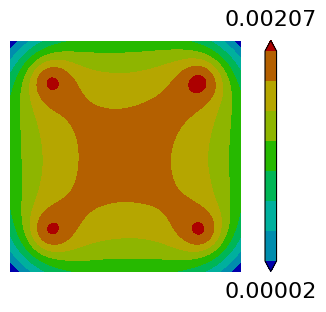} &
\hspace{-2pt}\includegraphics[height=\imgheight]{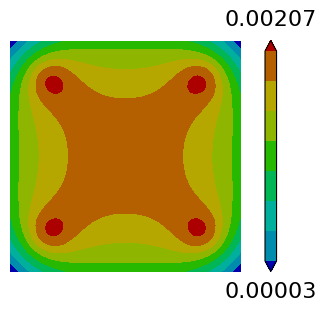} &
\hspace{-2pt}\includegraphics[height=\imgheight]{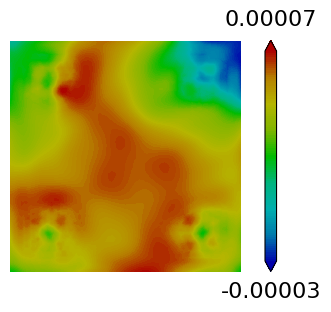} \\[-1pt]

\raisebox{0.9\imgheighth}{\textbf{Case 5}} &
\includegraphics[height=\imgheight]{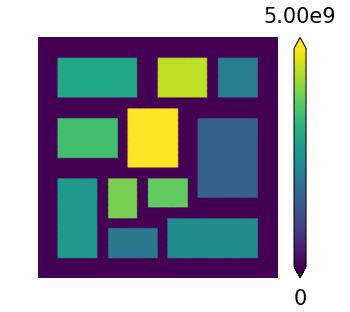} &
\hspace{-2pt}\includegraphics[height=\imgheight]{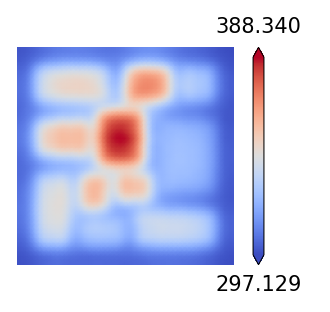} &
\hspace{-2pt}\includegraphics[height=\imgheight]{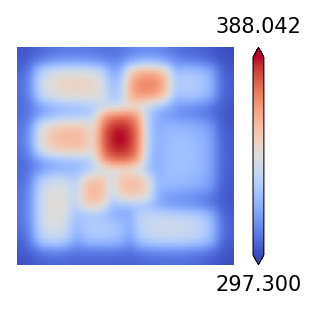} &
\hspace{-2pt}\includegraphics[height=\imgheight]{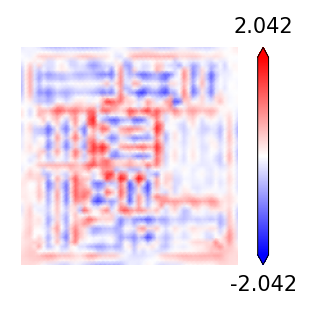} &
\hspace{-2pt}\includegraphics[height=\imgheight]{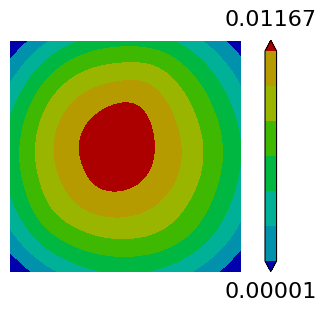} &
\hspace{-2pt}\includegraphics[height=\imgheight]{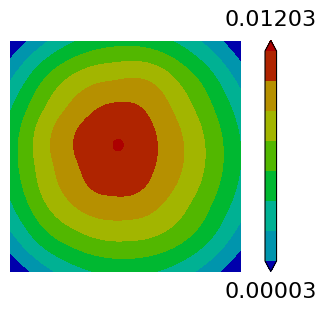} &
\hspace{-2pt}\includegraphics[height=\imgheight]{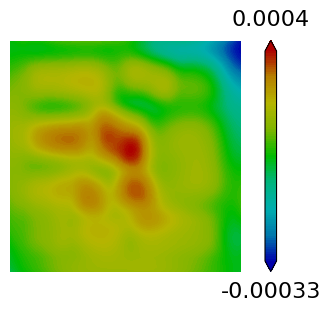} \\[-1pt]

\raisebox{0.9\imgheighth}{\textbf{Case 6}} &
\includegraphics[height=\imgheight]{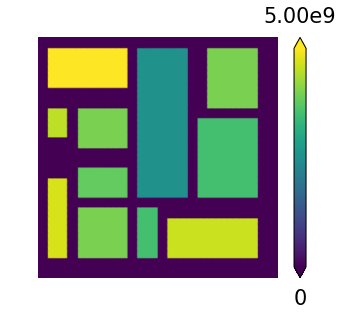} &
\hspace{-2pt}\includegraphics[height=\imgheight]{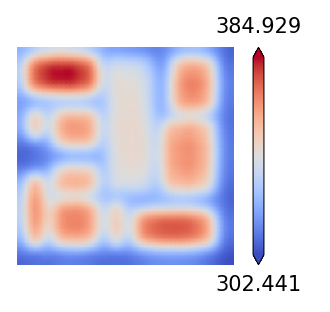} &
\hspace{-2pt}\includegraphics[height=\imgheight]{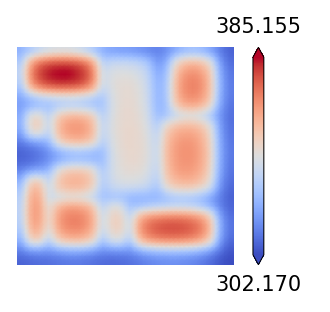} &
\hspace{-2pt}\includegraphics[height=\imgheight]{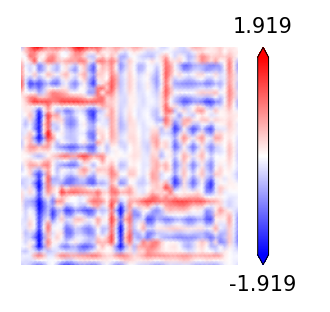} &
\hspace{-2pt}\includegraphics[height=\imgheight]{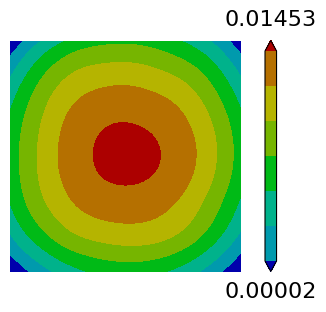} &
\hspace{-2pt}\includegraphics[height=\imgheight]{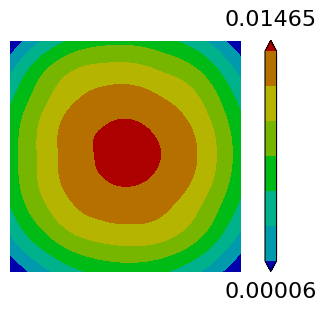} &
\hspace{-2pt}\includegraphics[height=\imgheight]{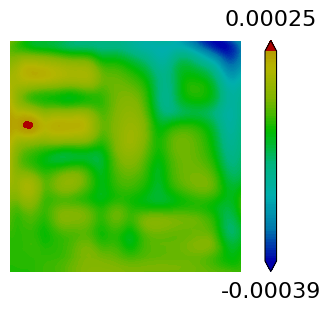} \\[-1pt]

\raisebox{0.9\imgheighth}{\textbf{Case 7}} &
\includegraphics[height=\imgheight]{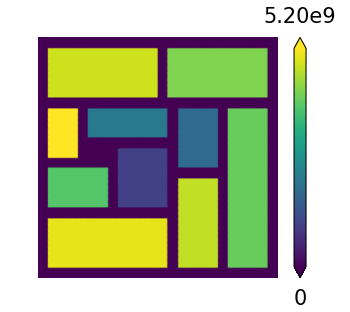} &
\hspace{-2pt}\includegraphics[height=\imgheight]{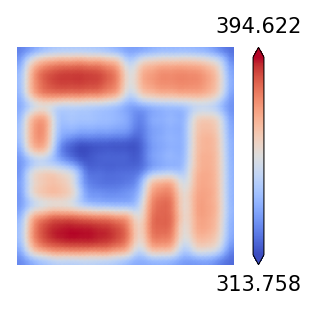} &
\hspace{-2pt}\includegraphics[height=\imgheight]{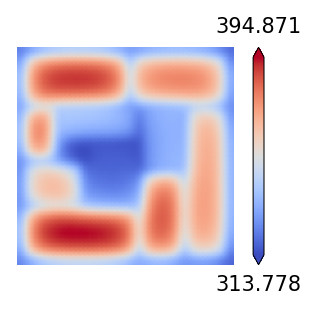} &
\hspace{-2pt}\includegraphics[height=\imgheight]{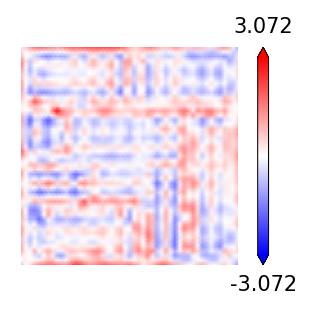} &
\hspace{-2pt}\includegraphics[height=\imgheight]{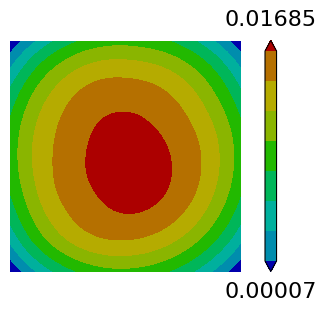} &
\hspace{-2pt}\includegraphics[height=\imgheight]{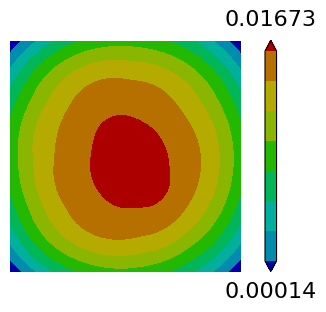} &
\hspace{-2pt}\includegraphics[height=\imgheight]{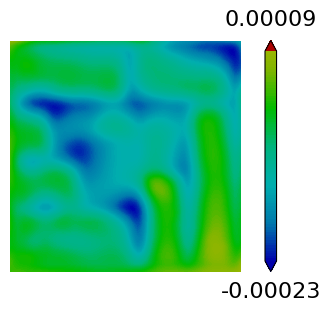} \\[-1pt]

\raisebox{0.9\imgheighth}{\textbf{Case 8}} &
\includegraphics[height=\imgheight]{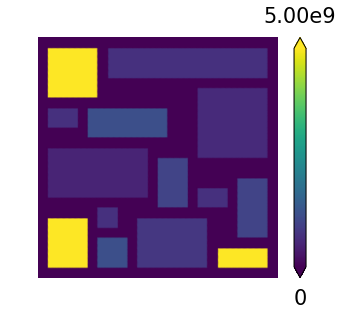} &
\hspace{-2pt}\includegraphics[height=\imgheight]{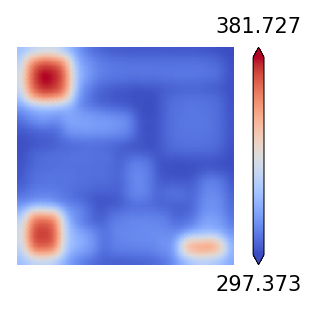} &
\hspace{-2pt}\includegraphics[height=\imgheight]{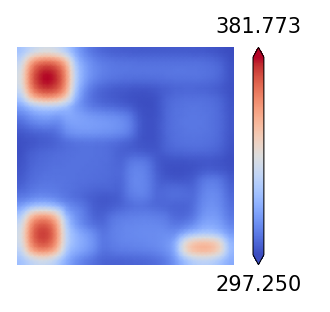} &
\hspace{-2pt}\includegraphics[height=\imgheight]{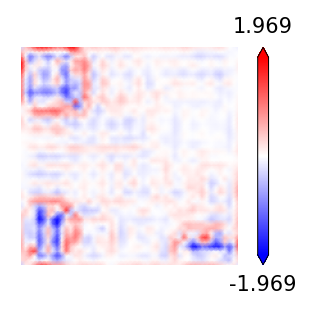} &
\hspace{-2pt}\includegraphics[height=\imgheight]{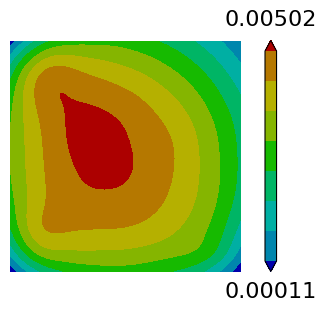} &
\hspace{-2pt}\includegraphics[height=\imgheight]{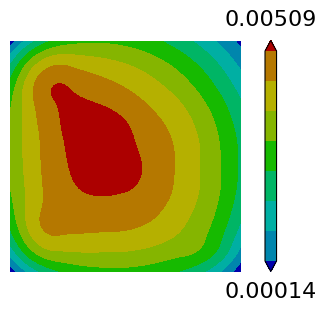} &
\hspace{-2pt}\includegraphics[height=\imgheight]{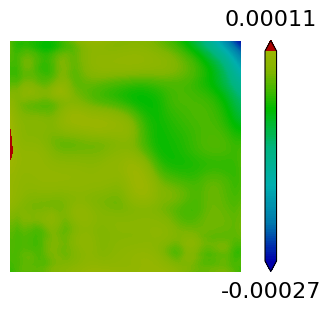} \\[2pt]

& \textbf{(a)} & \textbf{(b)} & \textbf{(c)} & \textbf{(d)} &
  \textbf{(e)} & \textbf{(f)} & \textbf{(g)} \\

\end{tabular}

\caption{Comparisons of temperature and warpage predictions for Case 1-8. (a) Power density. Temperature predicted by (b) T-PINN and (c) COMSOL. (d) Temperature error. Warpage predicted by (e) M-PINN and (f) COMSOL. (g) Warpage error.}
\label{fig:rscomp}
\end{figure*}

\begin{table}[htbp]
\centering
\caption{Material Properties for the Three-Layer Structure}
\label{tab:materials}
\begin{tabular}{lcccc}
\toprule
Layer & $z$-range (mm) & $E$ (GPa) & $\nu$ & $\alpha$ (ppm/K)  \\
\midrule
substrate & 0.0--0.8 & 400.0 & 0.14& 3.715 \\
interposer & 0.8--1.5 & 162.7 & 0.28 & 4 \\
chiplets & 1.5--2.25 & 129.6 & 0.28 & 4.21  \\
\bottomrule
\end{tabular}
\medskip
\end{table}

Fig. \ref{fig:3ds} illustrates a representative TSMC CoWoS-R advanced package comprising a GPU, a System-on-Chip (SoC), and six High Bandwidth Memory (HBM) modules. The 3D chiplet structure consists of three layers with overall dimensions of $24 \times 24 \times 2.25$ mm, where the layer interfaces are positioned at $z = 0.8$ mm and $z = 1.5$ mm. The material properties for the silicon-interposer-substrate assembly are summarized in Table \ref{tab:materials}. For the Stage-I T-PINN (temperature prediction), a uniform thermal conductivity of $60.0$ W/(m·K) is assumed. In contrast, the Stage-II M-PINN (warpage analysis) utilizes layer-dependent thermo-mechanical properties, including Young’s modulus ($E$), Poisson’s ratio ($\nu$), and the coefficient of thermal expansion ($\alpha$). To ensure numerical stability and continuity across interfaces, a hyperbolic tangent ($\tanh$) smoothing function is employed.

The T-PINN is implemented based on the ASRR-PINN architecture \cite{Zhou:DAC'25}. For the M-PINN, a uniform spatial point distribution is sampled across the 3D domain to capture the structural response, as depicted in Fig. \ref{fig:mesh}. The model is trained for 40,000 epochs using a weighted energy-based loss function optimized via the Adam algorithm. To enhance convergence robustness, a hybrid supervision strategy is adopted: the residual loss is monitored alongside the energy-based loss, and model parameters are updated only when both metrics show a concurrent decrease.

The prescribed power density distribution (Fig. \ref{fig:power}) reflects typical operational conditions, with higher power dissipation in the SoC region compared to the HBM modules. The steady-state temperature profile generated by the T-PINN (Fig. \ref{fig:temp}) subsequently serves as the non-uniform thermal load input for the Stage-II M-PINN to predict the resulting package warpage.

The predictive accuracy of the framework is validated along two critical cross-sectional paths. Fig. \ref{fig:line1} presents the diagonal displacement profile ($x = y$), demonstrating that the PINN predictions align closely with the COMSOL finite element reference. Similarly, Fig. \ref{fig:line2} indicates strong agreement along the vertical line at $x = 12$ mm. These results confirm the framework's capability to accurately capture complex 3D thermo-mechanical deformations in advanced packaging architectures with high fidelity.

\subsection{Accuracy and efficiency of the Proposed WarpagePINN}


\begin{table}[htbp]
\centering
\caption{Comparisons of WarpagePINN and COMSOL}
\label{tab:comparison}
\begin{tabular}{ccccc}
\toprule
Case & \begin{tabular}{@{}c@{}}T-MAE\\ (K)\end{tabular}  & \begin{tabular}{@{}c@{}}W-MAE\\ (mm)\end{tabular} & \begin{tabular}{@{}c@{}}WarpagePINN\\Inference Time (s)\end{tabular} & \begin{tabular}{@{}c@{}}COMSOL\\Time (s)\end{tabular}  \\
\midrule
1 & 0.1214 & 0.000039 & 0.00229 (8725$\times$) & 20 \\
2 & 0.0676 & 0.000046 & 0.00291 (8568$\times$) & 25 \\
3 & 0.1516 & 0.000030 & 0.00315 (7930$\times$) & 25 \\
4 & 0.1772 & 0.000042 & 0.00295 (7443$\times$) & 22 \\
5 & 0.1547 & 0.000160 & 0.00295 (8446$\times$) & 25 \\
6 & 0.1951 & 0.000225 & 0.00292 (8559$\times$) & 25 \\
7 & 0.2477 & 0.000079 & 0.00295 (8453$\times$) & 25 \\
8 & 0.0912 & 0.000057 & 0.00292 (8540$\times$) & 25 \\
\bottomrule
\end{tabular}
\medskip
\end{table}



To evaluate the accuracy and computational efficiency of the proposed WarpagePINN framework, a comprehensive performance analysis was conducted on eight distinct chiplet layouts with varying power density distributions, as illustrated in Fig.~\ref{fig:rscomp}. These cases maintain a consistent three-layer structural configuration to ensure a fair comparison. Table~\ref{tab:comparison} summarizes the quantitative metrics, including Mean Absolute Error (MAE) for temperature and warpage, as well as inference times. The results demonstrate that WarpagePINN achieves high-fidelity predictions across all scenarios. Specifically, in the worst-performing case (Case 7), the MAE for temperature and warpage remains as low as 0.25 K and 0.079 $\mu$m, respectively. In the best-case scenario (Case 2), the temperature error is further reduced to 0.07 K. Remarkably, while the conventional COMSOL Finite Element Analysis (FEA) requires approximately 20 to 25 s per simulation, WarpagePINN completes the inference in approximately 0.003 s. This achieves a computational speedup of up to 8725$\times$, effectively enabling real-time thermal-mechanical analysis.

The visual comparison of the temperature and warpage profiles at the chiplet layer for all eight layouts is presented in Fig.~\ref{fig:rscomp}. Column (a) illustrates the input power density distributions for the different layouts. Columns (b) and (e) represent the predicted temperature and warpage fields by WarpagePINN, which show an excellent agreement with the COMSOL ground truth results in columns (c) and (f). The error residuals shown in columns (d) and (g) further confirm the high visual and numerical fidelity of the proposed framework. The consistency between the predicted spatial gradients and the reference FEA results indicates that the physics-informed constraints effectively guide the model to capture the complex thermal-mechanical coupling effects without the need for dense mesh iterations.

Detailed morphological analysis of the warpage profiles reveals the sensitivity of the mechanical response to the heat source distribution. For the $2 \times 2$ chiplet layouts in Cases 1 and 2, the framework accurately captures the deformation shift. A single active chiplet (Case 1) results in maximum warpage concentrated at the chiplet center, whereas four active chiplets (Case 2) produce an asymmetric star-like profile due to thermal interaction. In the $5 \times 5$ uniform layouts, Case 3 exhibits an oval-shaped warpage contour driven by active center-line chiplets, while the symmetric corner activation in Case 4 results in a balanced, symmetric star-shaped profile. For the heterogeneous layouts (Cases 5--8), which involve chiplets of varying sizes, the model successfully captures diverse deformation morphologies. While Cases 5 through 7 generally manifest circular warpage patterns, Case 8 reveals a distinct ``water-drop" shaped profile, reflecting the highly non-uniform thermal load. These results demonstrate the robustness of WarpagePINN in generalizing across various topological configurations and power densities.

\subsection{Parameterized WarpagePINN and Uncertainty Quantification}
\begin{figure}[ht]
\centering
\subfigure[]{
\includegraphics[width=0.45\linewidth]{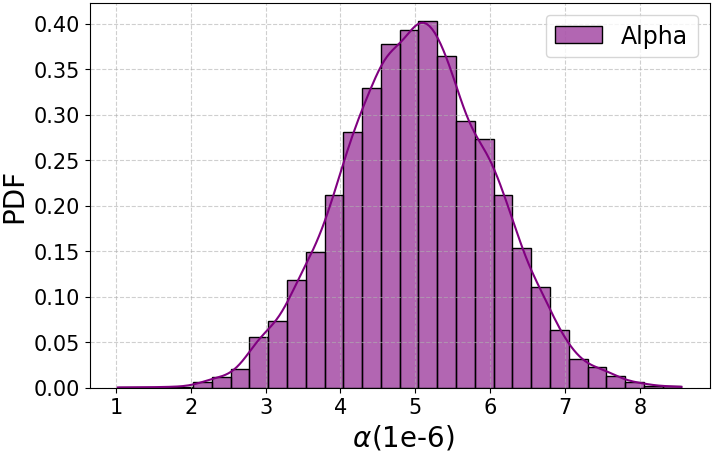}\label{fig:MCin}}
\subfigure[]{
\includegraphics[width=0.47\linewidth]
{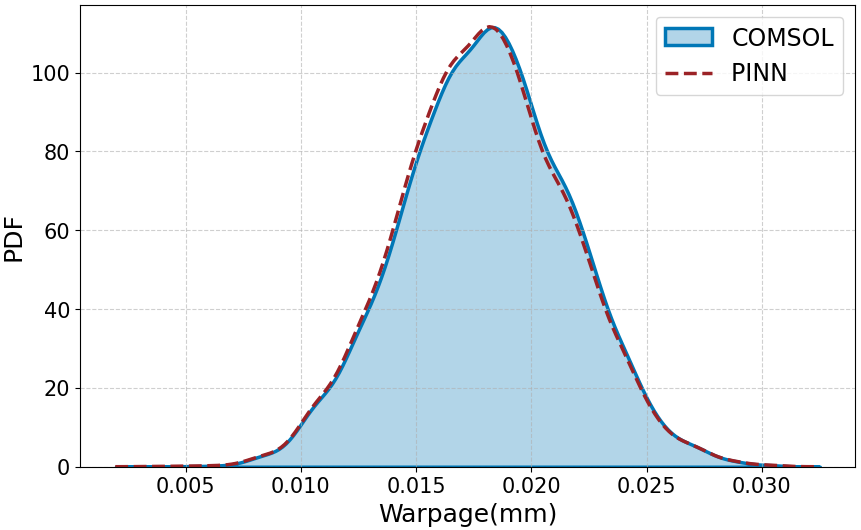}\label{fig:MCout}}
\caption{Monte Carlo simulation results with 10000 sampling points. (a) Distribution of the parameterized CTE inputs. (b) Corresponding predicted warpage outputs.}
\label{fig:wapage}
\end{figure}
The proposed framework is further extended into a parameterized WarpagePINN model, where the coefficient of thermal expansion (CTE), denoted as $\alpha$, is incorporated as an additional input dimension. This parameterization allows the network to learn a continuous mapping between material properties and mechanical responses. The model is initially trained using discrete CTE values spanning from 1 to 9 ppm/K. To evaluate the framework's generalizability and its capacity for uncertainty quantification, a Monte Carlo (MC) simulation was conducted. As illustrated in Fig. \ref{fig:MCin}, a total of 10,000 samples were generated from a normal distribution within the trained interval of 1--9 ppm/K to serve as inputs for the parametric inference.

A primary advantage of the parameterized WarpagePINN is its ability to perform near-instantaneous repeated inference once the offline training phase is complete. For the 10,000-sample MC analysis, the framework required a total cumulative inference time of only 59.8 s. In contrast, executing the equivalent number of high-fidelity simulations in COMSOL would require approximately 60,000 s (calculated based on average solver time per case). This represents a substantial computational speedup of approximately 1,000$\times$. Such efficiency gain is critical for design space exploration and sensitivity analysis, where traditional mesh-based solvers become computationally prohibitive.

In this analysis, warpage is defined as the maximum displacement differential across the chiplet domain. The statistical distribution of the predicted warpage results is presented in Fig. \ref{fig:MCout}. The output distribution predicted by the parameterized WarpagePINN demonstrates excellent agreement with the COMSOL reference results, accurately capturing the stochastic mechanical behavior induced by variations in CTE. These results confirm that the parameterized model not only maintains high numerical accuracy across a continuous parameter space but also provides the orders-of-magnitude acceleration necessary for rapid reliability assessment and iterative optimization in advanced packaging design.

\section{Conclusion}
\label{sec:concl}


This paper presents WarpagePINN, a novel physics-informed neural network framework for predicting thermal warpage in advanced packaging. The framework integrates two specialized components: T-PINN for temperature field prediction and M-PINN for mechanical warpage estimation. In the first stage, T-PINN employs a Fourier series representation that inherently satisfies boundary conditions, enabling training exclusively through a loss function derived from the governing heat equation. In the second stage, M-PINN utilizes a multilayer perceptron with a hybrid supervisory training strategy that optimizes an energy-based variational loss rather than the conventional PDE residual loss, significantly improving training stability and convergence. A parameterized extension of WarpagePINN is developed to quantify uncertainties associated with spatially varying coefficients of thermal expansion. Numerical experiments demonstrate that the proposed energy-based loss function achieves a 2.6$\times$ reduction in training time compared to the residual-based formulation. The parametric CTE study using WarpagePINN attains a 1000$\times$ computational speedup over finite-element analysis. Comprehensive validation across eight test cases shows that WarpagePINN matches FEM solutions with a mean absolute error of 0.2 $\mu$m, while delivering a remarkable 7400$\times$ speedup. These results establish WarpagePINN as an efficient and accurate surrogate for thermal-mechanical simulation in chiplet integration, offering substantial computational savings without compromising predictive fidelity.

\bibliographystyle{ieeetr}
\bibliography{./IEEErefer.bib}

\end{document}